\documentclass[english]{amsart}
\usepackage[T1]{fontenc}
\usepackage[cp850]{inputenc}
\usepackage{amsthm}
\usepackage{amssymb}

\numberwithin{equation}{section} %% Comment out for sequentially-numbered
\numberwithin{figure}{section} %% Comment out for sequentially-numbered

\theoremstyle{plain}
\newtheorem{thm}{Theorem}[section]
  \theoremstyle{definition}
  \newtheorem{defn}[thm]{Definition}
  \theoremstyle{remark}
  \newtheorem{claim}[thm]{Claim}
  \theoremstyle{remark}
  \newtheorem{rem}[thm]{Remark}
  \theoremstyle{plain}
  \newtheorem{prop}[thm]{Proposition}
  \theoremstyle{plain}
  \newtheorem{lem}[thm]{Lemma}
  \theoremstyle{plain}
  \newtheorem{cor}[thm]{Corollary}

\newtheorem{notation}[thm]{\sc Notation}

\allowdisplaybreaks

\usepackage{babel}

\begin{document}

\title{Quadratic Lie Algebras}

\author{Alessandro Ardizzoni}

\address{University of Ferrara, Department of Mathematics, Via Machiavelli
35, Ferrara, I-44100, Italy}

\email{alessandro.ardizzoni@unife.it}

\urladdr{http://www.unife.it/utenti/alessandro.ardizzoni}

\author{Fabio Stumbo}

\address{University of Ferrara, Department of Mathematics, Via Machiavelli
35, Ferrara, I-44100, Italy}

\email{f.stumbo@unife.it}

\urladdr{http://www.unife.it/utenti/fabio.stumbo}


\begin{abstract}
In this paper, the notion of universal enveloping algebra introduced
in [A. Ardizzoni, \emph{A First Sight Towards
Primitively Generated Connected Braided Bialgebras}, submitted. (arXiv:0805.3391v3)] is specialized to the case of braided vector
spaces whose Nichols algebra is quadratic as an algebra. In this setting
a classification of universal enveloping algebras for braided vector
spaces of dimension not greater than $2$ is handled. As an application, we investigate the structure of primitively generated connected braided bialgebras whose braided vector space of primitive elements forms a Nichols algebra which is quadratic algebra.
\end{abstract}


\thanks{This paper was written while the authors were members of GNSAGA with
partial financial support from MIUR (PRIN 2007)}

\maketitle
\tableofcontents{}

\section{Introduction}

Let $K$ be a fixed field. A braided vector space $(V,c)$ consists
of a vector space $V$ over $K$ and a $K$-linear map $c:V\otimes V\rightarrow V\otimes V$,
called braiding of $V$, which satisfies the quantum Yang-Baxter equation
\eqref{ec: braided equation}. A braided bialgebra is then a braided
vector space which is both an algebra and a coalgebra with structures
compatible with the braiding. Examples of braided bialgebras arises
as bialgebras in braided monoidal categories. 

In \cite{Ar- Universal}, the structure of primitively generated connected
braided bialgebras $A$ with respect to the braided vector space $\left(P,c\right)$
consisting of their primitive elements is investigated. The braided
vector space $\left(P,c\right)$ is endowed with a suitable $K$-linear
map $b:E\left(P,c\right)\rightarrow P$ and the datum $\left(P,c,b\right)$
is called a braided Lie algebra (ordinary Lie algebras
can be understood as a particular case of this notion). Here $E\left(P,c\right)$
denotes the space of primitive elements of the tensor algebra $T\left(P\right)$
having degree at least two. When the Nichols algebra of $\left(P,c\right)$
is obtained dividing out the tensor algebra $T(P)$ by the two-sided
ideal generated by $E\left(P,c\right)$, in \cite{Ar- Universal}
it is shown that $A$ can be recovered as a sort of universal enveloping
algebra $U\left(P,c,b\right)$ of the braided Lie algebra $\left(P,c,b\right)$.
As an application of this construction, in \cite{Ar- Universal} it
is proved that if $A$ is a connected braided bialgebra such that
the graded coalgebra $\mathrm{gr}A$ associated to the coradical filtration
is a quadratic algebra with respect to its natural braided bialgebra
structure, then, see Theorem \ref{teo: quadratic}, $A$ is isomorphic as
a braided bialgebra to the universal enveloping algebra $U\left(P,c,b\right)$
of the braided Lie algebra $\left(P,c,b\right)$ of primitive elements
of $A$. It is remarkable that in this context the Nichols algebra
$\emph{B}\left(P,c\right)$ is a quadratic algebra.

Motivated by this observation, in this paper we investigate the structure
of braided Lie algebras $\left(V,c,b\right)$ whose Nichols algebra
is quadratic as an algebra and of the corresponding universal enveloping
algebra $U\left(V,c,b\right)$. First, in Theorem \ref{thm:QLie},
we prove that when the Nichols algebra is quadratic, then the universal
enveloping algebra simplifies as follows: \[
U\left(V,c,b\right)=\frac{T\left(V,c\right)}{\left(\left(\mathrm{Id}-\bar{\beta}\right)\left[E_{2}\left(V,c\right)\right]\right)}\]
where $\bar{\beta}:E_{2}\left(V,c\right)\rightarrow V$ is the restriction
of $b$ to the space $E_{2}\left(V,c\right)=E\left(P,c\right)\cap V^{\otimes2}=\ker\left(c+\mathrm{Id}_{V^{\otimes2}}\right)$,
so that the datum $(V,c,\bar{\beta})$ completely encodes the structure
of the universal enveloping algebra and it is called a quadratic Lie
algebra (QLie algebra for short), see Definition \ref{def:QLie}.
Thus it is natural to write $U_{Q}\left(V,c,\bar{\beta}\right)$ instead
of $U\left(V,c,b\right)$. Now, in Lemma \ref{lem: lifting beta},
we prove that, under the further assumption that the braiding $c$
is root of a polynomial $f\in K\left[X\right]$ of the form $f=(X+1)h$ where $h\in K\left[X\right]$ is such that $h(-1)\neq 0$, then $\bar{\beta}:E_{2}\left(V,c\right)\rightarrow V$
is completely determined by the map $\beta:=\bar{\beta}h(c):V\otimes V\rightarrow V$.
Moreover, by Remark \ref{rem:UQ}, we have the following further simplification\[
U_{Q}\left(V,c,\bar{\beta}\right)=\frac{T\left(V,c\right)}{\left(h\left(c\right)\left(z\right)-\beta\left(z\right)\mid z\in V^{\otimes2}\right)}.\] Thus,
as above, the datum $(V,c,\beta)$ completely encodes the structure
of the universal enveloping algebra and it is called a lifted quadratic
Lie algebra (lifted QLie algebra for short), see Definition \ref{def:liftedQLie},
so we will write $U_{Q}\left(V,c,\beta\right)$ instead of $U_{Q}\left(V,c,\bar{\beta}\right)$.

In Section \ref{sec:Low-dimensional-cases}, we classify up to isomorphism
all lifted QLie algebras $\left(V,c,\beta\right)$ such that $\mathrm{dim}\left(V\right)\leq2$
and $\mathrm{dim}\left(\mathrm{Im}\beta\right)=1$, under mild assumptions
(the case $\mathrm{dim}\left(\mathrm{Im}\beta\right)=2$ is partially
treated in the Appendix \ref{sec:appendix}). The corresponding enveloping
algebras are also described explicitly in the main Theorem \ref{thm:teoClassif}.
As an application, in Theorem \ref{thm:sunto}, we study primitively
generated connected braided bialgebras $A$ with braided vector space
of primitive elements $(P,c)$ such that the Nichols algebra $\emph{B}\left(P,c\right)$
is a quadratic algebra and $c$ is root of a polynomial $f\in K\left[X\right]$
having $-1$ as a simple root. Such an $A$ is isomorphic to the universal
enveloping algebra $U_{Q}\left(V,c,\beta\right)$ of a suitable lifted
QLie algebra $\left(P,c,\beta\right)$ constructed on $\left(P,c\right)$.
Thus the previous results apply. 

We point out that many computations of this last part of the paper
have been handled with the help of the Computer Algebra System AXIOM
\cite{Axiom}.

\section{Preliminaries}

Throughout this paper $K$ will denote a field. All vector spaces
will be defined over $K$ and the tensor product of two vector spaces
will be denoted by $\otimes$.
\begin{thm}
Let $V$ be a vector space. Given a $K$-linear map $\alpha$ acting
on the space $V^{\otimes l},$ then we will denote by $\alpha_{i}$
then $K$-linear map $V^{\otimes i-1}\otimes\alpha\otimes V^{n-i-l+1}$
acting on the space $V^{\otimes n},n\geq l+i-1$.\end{thm}
\begin{defn}
Let $V$ be a vector space over a field $K$. A $K$-linear map $c:V\otimes V\rightarrow V\otimes V$
is called a \textbf{braiding of $V$} if it satisfies the quantum
Yang-Baxter equation\begin{equation}
c_{1}c_{2}c_{1}=c_{2}c_{1}c_{2}\label{ec: braided equation}\end{equation}
on $V\otimes V\otimes V$, where $c_{1}:=c\otimes V$ and $c_{2}:=V\otimes c.$
The pair $\left(V,c\right)$ will be called a\textbf{\ braided vector
space} (or $YB$-space). A morphism of braided vector spaces $(V,c_{V})$
and $(W,c_{W})$ is a $K$-linear map $f:V\rightarrow W$ such that
$c_{W}(f\otimes f)=(f\otimes f)c_{V}.$ 
\end{defn}
A general method for producing braided vector spaces is to take an
arbitrary braided category $(\mathcal{M},\otimes,K,a,l,r,c),$ which
is a monoidal subcategory of the category of $K$-vector spaces. Hence
any object $V\in\mathcal{M}$ can be regarded as a braided vector
space with respect to $c:=c_{V,V}.$ Here, $c_{X,Y}:X\otimes Y\rightarrow Y\otimes X$
denotes the braiding in $\mathcal{M}.$ The category of comodules
over a coquasitriangular Hopf algebra and the category of Yetter-Drinfeld
modules are examples of such categories. More particularly, every
bicharacter of a group $G$ induces a braiding on the category of
$G$-graded vector spaces.
\begin{defn}
\cite{Ba} The quadruple $(A,m_{A},u_{A},c_{A})$ is called a \textbf{braided
algebra} if
\begin{itemize}
\item $(A,m_{A},u_{A})$ is an associative unital algebra;
\item $(A,c_{A})$ is a braided vector space;
\item $m_{A}$ and $u_{A}$ commute with $c_{A}$, that is the following
conditions hold: \begin{gather}
c_{A}(m_{A}\otimes A)=(A\otimes m_{A})(c_{A}\otimes A)(A\otimes c_{A}),\label{Br2}\\
c_{A}(A\otimes m_{A})=(m_{A}\otimes A)\left(A\otimes c_{A}\right)(c_{A}\otimes A),\label{Br3}\\
c_{A}(u_{A}\otimes A)=A\otimes u_{A},\qquad c_{A}(A\otimes u_{A})=u_{A}\otimes A.\label{Br4}\end{gather}

\end{itemize}
A morphism of braided algebras is, by definition, a morphism of ordinary
algebras which, in addition, is a morphism of braided vector spaces.
\end{defn}
\begin{defn}
The quadruple $(C,\Delta_{C},\varepsilon_{C},c_{C})$ is called a
\textbf{braided coalgebra} if
\begin{itemize}
\item $(C,\Delta_{C},\varepsilon_{C})$ is a coassociative counital coalgebra;
\item $(C,c_{C})$ is a braided vector space;
\item $\Delta_{C}$ and $\varepsilon_{C}$ commute with $c_{C}$, that is
the following relations hold: \begin{gather}
(\Delta_{C}\otimes C)c_{C}=(C\otimes c_{C})(c_{C}\otimes C)(C\otimes\Delta_{C}),\label{Br5}\\
(C\otimes\Delta_{C})c_{C}=(c_{C}\otimes C)(C\otimes c_{C})(\Delta_{C}\otimes C),\label{Br6}\\
(\varepsilon_{C}\otimes C)c_{C}=C\otimes\varepsilon_{C},\qquad(C\otimes\varepsilon_{C})c_{C}=\varepsilon_{C}\otimes C.\label{Br7}\end{gather}

\end{itemize}
A morphism of braided coalgebras is, by definition, a morphism of
ordinary coalgebras which, in addition, is a morphism of braided vector
spaces.
\end{defn}
\begin{defn}
\cite[Definition 5.1]{Ta} A sextuple $(B,m_{B},u_{B},\Delta_{B},\varepsilon_{B},c_{B})$
is a called a \textbf{braided bialgebra} if
\begin{itemize}
\item $(B,m_{B},u_{B},c_{B})$ is a braided algebra
\item $(B,\Delta_{B},\varepsilon_{B},c_{B})$ is a braided coalgebra
\item the following relations hold:\begin{equation}
\Delta_{B}m_{B}=(m_{B}\otimes m_{B})(B\otimes c_{B}\otimes B)(\Delta_{B}\otimes\Delta_{B}).\label{Br1}\end{equation}

\end{itemize}
\end{defn}
\begin{defn}
\label{def: graded bialg}We will need graded versions of braided
algebras, coalgebras and bialgebras. A \textbf{graded braided algebra}
is a braided algebra $(A,m_{A},u_{A},c_{A})$ such that $A=\bigoplus_{n\in\mathbb{N}}A^{n}$
and $m_{A}(A^{m}\otimes A^{n})\subseteq A^{m+n},$ for every $m,n\in\mathbb{N}$.
The braiding $c_{A}$ is assumed to satisfy $c_{A}(A^{m}\otimes A^{n})\subseteq A^{n}\otimes A^{m}.$
It is easy to see that $1_{A}=u_{A}\left(1_{K}\right)\in A^{0}$.
Therefore a graded braided algebra is defined by maps $m_{A}^{m,n}:A^{m}\otimes A^{n}\rightarrow A^{m+n}$
and $c_{A}^{m,n}:A^{m}\otimes A^{n}\rightarrow A^{n}\otimes A^{m}$,
and by an element $1\in A^{0}$ such that for all $n,m,p\in\mathbb{N},a\in A^{n}$
\begin{gather}
m_{A}^{n+m,p}(m_{A}^{n,m}\otimes A^{p})=m_{A}^{n,m+p}(A^{n}\otimes m_{A}^{m,p}),\label{gr1}\\
m_{A}^{0,n}(1\otimes a)=a=m_{A}^{n,0}(a\otimes1),\label{gr2}\\
c_{A}^{n+m,p}(m_{A}^{n,m}\otimes A^{p})=(A^{p}\otimes m_{A}^{n,m})(c_{A}^{n,p}\otimes A^{m})(A^{n}\otimes c_{A}^{m,p}),\label{gr3}\\
c_{A}^{n,m+p}(A^{n}\otimes m_{A}^{m,p})=(m_{A}^{m,p}\otimes A^{n})(A^{m}\otimes c_{A}^{n,p})(c_{A}^{n,m}\otimes A^{p}),\label{gr4}\\
c_{A}^{0,n}(1\otimes a)=a\otimes1\qquad\text{and\qquad}\mathfrak{\ }c_{A}^{n,0}(a\otimes1)=1\otimes a.\label{gr5}\end{gather}
The multiplication $m_{A}$ can be recovered from $(m_{A}^{n,m})_{n,m\in\mathbb{N}}$
as the unique $K$-linear map such that $m_{A}(x\otimes y)=m_{A}^{p,q}(x\otimes y),\text{ for all }p,q\in\mathbb{N},x\in A^{p},y\in A^{q}.$
Analogously, the braiding $c_{A}$ is uniquely defined by $c_{A}(x\otimes y)=c_{A}^{p,q}(x\otimes y),\text{ for all }p,q\in\mathbb{N},x\in A^{p},y\in A^{q}.$
We will say that $m_{A}^{n,m}$ and $c_{A}^{n,m}$ are the $(n,m)$-homogeneous
components of $\nabla$ and $c_{A}$, respectively.

\textbf{Graded braided coalgebras} can by described in a similar way.
By definition a braided coalgebra $(C,\Delta_{C},\varepsilon_{C},c_{C})$
is graded if $C=\bigoplus_{n\in\mathbb{N}}C^{n},$ $\Delta_{C}(C^{n})\subseteq\sum_{p+q=n}C^{p}\otimes C^{q}$,
$c_{C}(C^{n}\otimes C^{m})\subseteq C^{m}\otimes C^{n}$ and $\varepsilon_{C|_{C_{n}}}=0$,
for $n>0$ . If $\pi^{p}$ denotes the projection onto $C^{p}$ then
the comultiplication $\Delta_{C}$ is uniquely defined by maps $\Delta_{C}^{p,q}:C^{p+q}\rightarrow C^{p}\otimes C^{q}$,
where $\Delta_{C}^{p,q}:=(\pi^{p}\otimes\pi^{q})\Delta_{C}|_{C^{p+q}}$.
The counit is given by a map $\varepsilon_{C}^{0}:C^{0}\rightarrow K,$
while the braiding $c_{C}$ is uniquely determined by a family $(c_{C}^{n,m})_{n,m\in\mathbb{N}},$
as for braided algebras. The families $(\Delta_{C}^{n,m})_{n,m\in\mathbb{N}},$
$(c_{C}^{n,m})_{n,m\in\mathbb{N}}$ and $\varepsilon_{C}^{0}$ has
to satisfy the relations that are dual to (\ref{gr1}) -- (\ref{gr5}),
namely for all $n,m,p\in\mathbb{N},c\in C^{n},d\in C^{0}$: \begin{gather}
(\Delta_{C}^{n,m}\otimes C^{p})\Delta_{C}^{n+m,p}=(C^{n}\otimes\Delta_{C}^{m,p})\Delta_{C}^{n,m+p},\label{c1}\\
(\varepsilon_{C}^{0}\otimes C^{n})\Delta_{C}^{0,n}(c)=c=(C^{n}\otimes\varepsilon_{C}^{0})\Delta_{C}^{n,0}(c),\label{c2}\\
(C^{p}\otimes\Delta_{C}^{n,m})c_{C}^{n+m,p}=(c_{C}^{n,p}\otimes C^{m})(C^{n}\otimes c_{C}^{m,p})(\Delta_{C}^{n,m}\otimes C^{p}),\label{c3}\\
(\Delta_{C}^{m,p}\otimes C^{n})c_{C}^{n,m+p}=(C^{m}\otimes c_{C}^{n,p})(c_{C}^{n,m}\otimes C^{p})(C^{n}\otimes\Delta_{C}^{m,p}),\label{c4}\\
(\varepsilon_{C}^{0}\otimes C)c_{C}(c\otimes d)=\varepsilon_{C}^{0}(d)c=(C\otimes\varepsilon_{C}^{0})c_{C}(d\otimes c).\label{c5}\end{gather}
We will say that $\Delta_{C}^{n,m}$ is the $(n,m)$-homogeneous component
of $\Delta_{C}$.

A \textbf{graded braided bialgebra} is a braided bialgebra which is
graded both as an algebra and as a coalgebra. 
\end{defn}

\begin{defn}
Recall that a coalgebra $C$ is called \textbf{connected} if $C_{0}$,
the coradical of $C$ (i.e the sum of simple subcoalgebras of $C$),
is one dimensional. In this case there is a unique group-like element
$1_{C}\in C$ such that $C_{0}=K1_{C}$. A morphism of connected coalgebras
is a coalgebra homomorphisms (clearly it preserves the grouplike elements). 

By definition, a braided coalgebra $\left(C,c_{C}\right)$ is \textbf{connected}
if $C_{0}=K1_{C}$ and, for any $x\in C$, \begin{equation}
c_{C}(x\otimes1_{C})=1_{C}\otimes x,\qquad c_{C}(1_{C}\otimes x)=x\otimes1_{C}.\label{de: connected}\end{equation}
\end{defn}

\subsection{The universal enveloping algebra}

In this subsection we recall some definitions and notations from \cite{Ar- Universal}.
\begin{defn}
\label{def: E(C)}Let $T:=T\left(V,c\right)$ and set \[
E\left(V,c\right):=\bigoplus\limits _{n\in\mathbb{N}}E_{n}\left(V,c\right)\]
where\[
E_{n}\left(V,c\right):=\left\{ \begin{tabular}{ll}
 $0$  &  if $n=0,1,$ \\
$\bigcap\limits _{\substack{a,b\geq1\\
a+b=n}
}\ker\left(\Delta_{T}^{a,b}\right)$  &  if $n\geq2.$\end{tabular}\right.\]
The elements of $E\left(V,c\right)$ are the primitive elements in
$T$ of degree greater then $1$ (cf. \cite[Lemma 3.3]{Ar- Universal}).
Note that $E_{2}\left(V,c\right)=\ker\left(\Delta_{T}^{1,1}\right)=\ker\left(c+\mathrm{Id}_{V^{\otimes2}}\right)$.\end{defn}
\begin{claim}
It is easy to check that the braiding $c_{T}$ of the tensor algebra
$T:=T\left(V,c\right)$ induces for all $u\in\mathbb{N}$ maps\begin{eqnarray*}
c_{E_{t}\left(V,c\right),V^{\otimes u}}:E_{t}\left(V,c\right)\otimes V^{\otimes u} & \rightarrow & V^{\otimes u}\otimes E_{t}\left(V,c\right),\\
c_{V^{\otimes u},E_{t}\left(V,c\right)}:V^{\otimes u}\otimes E_{t}\left(V,c\right) & \rightarrow & E_{t}\left(V,c\right)\otimes V^{\otimes u}.\end{eqnarray*}
Thus we can define the following maps\begin{eqnarray*}
c_{E\left(V,c\right),V^{\otimes u}} & : & =\bigoplus_{t\in\mathbb{N}}c_{E_{t}\left(V,c\right),V^{\otimes u}}:E\left(V,c\right)\otimes V^{\otimes u}\rightarrow V^{\otimes u}\otimes E\left(V,c\right),\\
c_{V^{\otimes u},E\left(V,c\right)} & : & =\bigoplus_{t\in\mathbb{N}}c_{V^{\otimes u},E_{t}\left(V,c\right)}:V^{\otimes u}\otimes E\left(V,c\right)\rightarrow E\left(V,c\right)\otimes V^{\otimes u}.\end{eqnarray*}
\end{claim}
\begin{defn}
A \textbf{(braided) bracket} on a braided vector space $\left(V,c\right)$
is a $K$-linear map $b:E\left(V,c\right)\rightarrow V$ such that
\begin{equation}
c\left(b\otimes V\right)=\left(V\otimes b\right)c_{E\left(V,c\right),V}\qquad c\left(V\otimes b\right)=\left(b\otimes V\right)c_{V,E\left(V,c\right)}.\label{form: c-bracket}\end{equation}
The restriction of $b$ to $E_{t}\left(V,c\right)$ will be denoted
by $b^{t}:E_{t}\left(V,c\right)\rightarrow V.$

If $b$ is a bracket on a braided vector space $\left(V,c\right),$
then we define \textbf{the universal enveloping algebra of} $\left(V,c,b\right)$
to be \begin{equation}
U\left(V,c,b\right):=\frac{T\left(V,c\right)}{\left(\left(\mathrm{Id}-b\right)\left[E\left(V,c\right)\right]\right)}\label{form: U(V,c,b)}\end{equation}
which, by \cite[Theorem 6.3]{AMS-MM}, is a braided bialgebra quotient
of the tensor algebra $T\left(V,c\right).$ Note that $0$ is always
a bracket on $\left(V,c\right)$ so that it makes sense to define
the \textbf{symmetric algebra of} $\left(V,c\right)$ by setting \[
S\left(V,c\right):=U\left(V,c,0\right)=\frac{T\left(V,c\right)}{\left(E\left(V,c\right)\right)}.\]
We will denote by $\pi_{U}:T\left(V,c\right)\rightarrow U\left(V,c,b\right)$
(resp. $\pi_{S}:T\left(V,c\right)\rightarrow S\left(V,c\right)$)
the canonical projection and by $i_{U}:V\rightarrow U\left(V,c,b\right)$
(resp. $i_{S}:V\rightarrow S\left(V,c\right)$) its restriction to
$V$. 

\end{defn}
\begin{defn}
\cite[Definition 5.6]{Ar- Universal}\label{def: c-Lie alg} We say
that $\left(V,c,b\right)$ is a \textbf{braided Lie algebra} whenever
\begin{itemize}
\item $\left(V,c\right)$ is a braided vector space;
\item $b:E\left(V,c\right)\rightarrow V$ is a bracket on $\left(V,c\right)$;
\item the canonical $K$-linear map $i_{U}:V\rightarrow U\left(V,c,b\right)$
is injective i.e. \[
V\cap\ker\left(\pi_{U}\right)=\ker\left(i_{U}\right)=0.\]
 
\end{itemize}
\end{defn}
\begin{claim}
\label{claim: Gamma}Let $\left(V,c\right)$ be a braided vector space
and let $T=T\left(V,c\right)$. By the universal property of the tensor
algebra there is a unique algebra homomorphism \[
\mathrm{\Gamma}^{T}:T\left(V,c\right)\rightarrow T^{c}\left(V,c\right)\]
such that $\mathrm{\Gamma}_{\mid V}^{T}=\mathrm{Id}_{V},$ where $T^{c}\left(V,c\right)$
denotes the quantum shuffle algebra. This is a morphism of graded
braided bialgebras. The \textbf{Nichols algebra} of $(V,c)$ is defined by\[
\emph{B}\left(V,c\right)=\text{Im}\left(\mathrm{\Gamma}^{T}\right)\simeq\frac{T\left(V,c\right)}{\ker\left(\mathrm{\Gamma}^{T}\right)}.\]

\end{claim}

\section{General results on quadratic algebras}

We recall from \cite[page 19]{Manin} the definition of quadratic
algebra.
\begin{defn}
\label{def: quadratic}A \textbf{quadratic algebra} is an associative
graded $K$-algebra $A=\oplus_{n\in\mathbb{N}}A^{n}$ such that:
\begin{enumerate}
\item [1)] $A^{0}=K$;
\item [2)] $A$ is generated as a $K$-algebra by $A^{1}$;
\item [3)] the ideal of relations among elements of $A^{1}$ is generated
by the subspace of all quadratic relations $R(A)\subseteq A^{1}\otimes A^{1}.$ 
\end{enumerate}
Equivalently $A$ is a graded $K$-algebra such that the natural map
$\pi:T_{K}(A^{1})\rightarrow A$ from the tensor algebra generated
by $A^{1}$ is surjective and $\mathrm{\ker}(\pi)$ is generated as
a two sided ideal in $T_{K}(A^{1})$ by $\mathrm{\ker}(\pi)\cap[A^{1}\otimes A^{1}].$ \end{defn}
\begin{thm}
\label{teo: quadratic}\cite[Theorem 9.5]{Ar- Universal} Let $\left(A,c_{A}\right)$
be a connected braided bialgebra such that the graded coalgebra associated
to the coradical filtration is a quadratic algebra with respect to
its natural braided bialgebra structure. Let $\left(P,c_{P},b_{P}\right)$
be the infinitesimal braided Lie algebra of $A$. Then $\emph{B}\left(P,c_{P}\right)$
is a quadratic algebra and $A$ is isomorphic to $U\left(P,c_{P},b_{P}\right)$
as a braided bialgebra. 
\end{thm}
Theorem \ref{teo: quadratic} justifies our interest in the study
of the universal enveloping algebra of braided Lie algebras $\left(V,c,b\right)$
such that $\emph{B}\left(V,c\right)$ is a quadratic algebra. Examples
of braided vector spaces $\left(V,c\right)$ such that $\emph{B}\left(V,c\right)$
is a quadratic algebra can be found e.g. in \cite{MS-Pointed indec}
and in \cite{AG- Racks}. By \cite[Proposition 3.4]{AS}, another
example is given by braided vector spaces of Hecke-type with regular
mark.

\subsection{Quadratic Lie algebras}

In this subsection we introduce and study the notion of quadratic
Lie algebra: we will see that this notion completely determines the
universal enveloping algebra of a braided Lie algebra.
\begin{defn}
A \textbf{quadratic bracket} (Qbracket for short) on a braided vector
space $\left(V,c\right)$ is a $K$-linear map $\overline{\beta}=\overline{\beta}_{V}:E_{2}\left(V,c\right)\rightarrow V$
such that \begin{equation}
c\overline{\beta}_{1}=\overline{\beta}_{2}c_{E_{2}\left(V,c\right),V}\qquad c\overline{\beta}_{2}=\overline{\beta}_{1}c_{V,E_{2}\left(V,c\right)}.\label{form:bracket}\end{equation}
A \textbf{morphism of Qbrackets} is a morphism of braided vector spaces
$f:\left(V,c_{V}\right)\rightarrow\left(W,c_{W}\right)$ such that
$f\circ\overline{\beta}_{V}=\overline{\beta}_{W}\circ\left(f\otimes f\right).$
For any braided vector space $\left(V,c\right)$, we set \[
\overline{E_{2}\left(V,c\right)}=\left(E_{2}\left(V,c\right)\otimes V\right)\cap\left(V\otimes E_{2}\left(V,c\right)\right)=\left\{ z\in V^{\otimes3}\mid c_{1}\left(z\right)=-z=c_{2}\left(z\right)\right\} .\]
\end{defn}
\begin{rem}
Let $\overline{\beta}:E_{2}\left(V,c\right)\rightarrow V$ be a Qbracket
on a braided vector space $\left(V,c\right).$ Then\begin{equation}
\left(\overline{\beta}_{1}-\overline{\beta}_{2}\right)\left(\overline{E_{2}\left(V,c\right)}\right)\subseteq E_{2}\left(V,c\right).\label{form:correctness}\end{equation}
In fact, let $z\in\overline{E_{2}\left(V,c\right)}.$ Then \begin{eqnarray*}
c\left(\overline{\beta}_{1}-\overline{\beta}_{2}\right)\left(z\right) & = & c\overline{\beta}_{1}\left(z\right)-c\overline{\beta}_{2}\left(z\right)\\
 & \overset{\text{(\ref{form:bracket})}}{=} & \overline{\beta}_{2}c_{E_{2}\left(V,c\right),V}\left(z\right)-\overline{\beta}_{1}c_{V,E_{2}\left(V,c\right)}\left(z\right)\\
 & = & \overline{\beta}_{2}\left(z\right)-\overline{\beta}_{1}\left(z\right).\end{eqnarray*}
 \end{rem}
\begin{defn}
\label{def:QLie}A \textbf{quadratic Lie algebra} (QLie algebra for
short) is a tern $\left(V,c,\overline{\beta}\right)$ where
\begin{itemize}
\item $\left(V,c\right)$ is a braided vector space;
\item $\overline{\beta}:E_{2}\left(V,c\right)\rightarrow V$ is a Qbracket
such that \begin{equation}
\overline{\beta}\left(\overline{\beta}_{1}-\overline{\beta}_{2}\right)\left(\overline{E_{2}\left(V,c\right)}\right)=0.\label{form:Jacobi}\end{equation}
 Note that \eqref{form:Jacobi} makes sense in view of \eqref{form:correctness}. 
\end{itemize}
A \textbf{morphism of QLie algebras} $f:\left(V,c_{V},\overline{\beta}_{V}\right)\rightarrow\left(W,c_{W},\overline{\beta}_{W}\right)$
is by definition a morphism of Qbrackets $f:\left(V,c_{V}\right)\rightarrow\left(W,c_{W}\right)$. 

If $\overline{\beta}$ is a Qbracket on a braided vector space $\left(V,c\right),$
then we define \textbf{the quadratic universal enveloping algebra
of} $\left(V,c,\overline{\beta}\right)$ to be \[
U_{Q}\left(V,c,\overline{\beta}\right):=\frac{T\left(V,c\right)}{\left(\left(\mathrm{Id}-\overline{\beta}\right)\left[E_{2}\left(V,c\right)\right]\right)}.\]
Mimicking the proof of \cite[Theorem 3.9]{Ar- Universal}, one easily
verifies that $U_{Q}\left(V,c,\overline{\beta}\right)$ is indeed
a braided bialgebra quotient of the tensor algebra $T\left(V,c\right).$
Note that $0$ is always a Qbracket on $\left(V,c\right)$ so that
it makes sense to define the \textbf{quadratic} \textbf{symmetric
algebra of} $\left(V,c\right)$ by setting \[
S_{Q}\left(V,c\right):=U_{Q}\left(V,c,0\right)=\frac{T\left(V,c\right)}{\left(E_{2}\left(V,c\right)\right)}.\]
We will denote by $\pi_{U_{Q}}:T\left(V,c\right)\rightarrow U_{Q}\left(V,c,b\right)$
(resp. $\pi_{S_{Q}}:T\left(V,c\right)\rightarrow S_{Q}\left(V,c\right)$)
the canonical projection and by $i_{U_{Q}}:V\rightarrow U_{Q}\left(V,c,b\right)$
(resp. $i_{S_{Q}}:V\rightarrow S_{Q}\left(V,c\right)$) its restriction
to $V$. \end{defn}

\begin{rem}
Let $(V,I\oplus I^{*}=V^{\otimes2},\overline{\beta})$ be a ``braided
Lie algebra'' in the sense of \cite[Definition 1]{Gurevich-Hecke}.
Then $\overline{\beta}$ is a map from $V\otimes V$ to $V$ where $V$ is an
object in the braided monoidal category $\mathfrak{A}$ defined therein.
Hence the braiding of $\mathfrak{A}$ induces a braiding
on $V$ that we denote by $c$. Suppose
$I:=E_{2}\left(V,c\right)=\mathrm{ker}(c+\mathrm{Id}_{V^{\otimes2}})$
(for instance this is what happens in \cite[bottom of page
325]{Gurevich-Hecke}
where $V'$ plays the role of our $V$, see also \cite[Example 9.6]{Ar-
Universal})
Then condition 0 in Gurevich's definition means that $S_{Q}\left(V,c\right)$
is a Koszul algebra. Condition 1 yields that $\overline{\beta}$ is completely
determined by its restriction to $I$. Condition 2 entails
that \eqref{form:Jacobi} holds. Condition 3 tells that $\overline{\beta}$
is a morphism in $\mathfrak{A}$ so that $\overline{\beta}$ is compatible
with the braiding of $\mathfrak{A}$ whence it is a Qbracket. Summing
up we get that $\left(V,c,\overline{\beta}\right)$ is a QLie algebra
such that $S_{Q}\left(V,c\right)$ is a Koszul algebra. Note also
that $U_{Q}\left(V,c,\overline{\beta}\right)=U\left(\mathfrak{g}\right)$
in Gurevich's sense.
\end{rem}

\begin{claim}
\label{claim: graded}Let $\left(V,c\right)$ be a braided vector
space and let $\overline{\beta}:E_{2}\left(V,c\right)\rightarrow V$
be a Qbracket. Set\[
T:=T\left(V,c\right)\qquad\text{and}\qquad U_{Q}:=U_{Q}\left(V,c,b\right).\]
Let $\pi_{U_{Q}}:T\rightarrow U_{Q}$ be the canonical projection.
By construction $\pi_{U_{Q}}$ is a morphism of braided bialgebras.
Set $T^{\leq n}:=\oplus_{0\leq t\leq n}V^{\otimes t}$ and \[
U_{n}^{\prime}:=\pi_{U_{Q}}\left(T^{\leq n}\right).\]
Then $\left(U_{n}^{\prime}\right)_{n\in\mathbb{N}}$ is both an algebra
and a coalgebra filtration on $U_{Q}$, which is called the \emph{standard
filtration on} $U_{Q}$. Note that this filtration is not the coradical
filtration of $U_{Q}$ in general. Still one has $U_{n}^{\prime}=\pi_{U_{Q}}\left(T^{\leq n}\right)\subseteq\pi_{U_{Q}}\left(T_{n}\right)\subseteq\left(U_{Q}\right)_{n}$
where $T_{n}$ and $\left(U_{Q}\right)_{n}$ denote the $n$-th terms
of the coradical filtration of $T$ and $U_{Q}$ respectively. Denote
by\[
\mathrm{gr}^{\prime}\left(U_{Q}\right):=\oplus_{n\in\mathbb{N}}\frac{U_{n}^{\prime}}{U_{n-1}^{\prime}}.\]
the graded coalgebra associated to the standard filtration, see \cite[page 228]{Sw}.

If $\overline{\beta}=0$, then $S_{Q}\left(V,c\right)=U_{Q}\left(V,c,0\right)$
is a graded bialgebra of the form $S_{Q}\left(V,c\right)=\oplus_{n\in\mathbb{N}}S_{Q}^{n}\left(V,c\right)$.
The standard filtration on $S_{Q}\left(V,c\right)$ is the filtration
associated to this grading. 
\end{claim}

\begin{prop}
\label{pro: theta}Let $\overline{\beta}:E_{2}\left(V,c\right)\rightarrow V$
be a Qbracket on a braided vector space $\left(V,c\right)$. Then
$\mathrm{gr}^{\prime}\left(U_{Q}\left(V,c,b\right)\right)$ is a graded
braided bialgebra and there is a canonical morphism of graded braided
bialgebras $\theta:S_{Q}\left(V,c\right)\rightarrow\mathrm{gr}^{\prime}\left(U_{Q}\left(V,c,b\right)\right)$
which is surjective and lifts the map $\theta_{1}:V\rightarrow U_{1}^{\prime}/U_{0}^{\prime}:v\mapsto\pi_{U_{Q}}\left(v\right)+U_{0}^{\prime}.$ \end{prop}
\begin{proof}
It is similar to the proof of \cite[Proposition 5.14]{Ar- Universal}. \end{proof}
\begin{thm}
\label{teo: magnum}Let $\overline{\beta}:E_{2}\left(V,c\right)\rightarrow V$
be a Qbracket on a braided vector space $\left(V,c\right)$. Then
$\emph{B}\left(V,c\right)$ is a quadratic algebra if and only if
$S_{Q}\left(V,c\right)=\emph{B}\left(V,c\right)$. If this holds,
the following assertions are equivalent.
\begin{enumerate}
\item [(i)] $i_{U_{Q}}:V\rightarrow U_{Q}\left(V,c,\overline{\beta}\right)$
is injective.
\item [(ii)] $V\cap\ker\left(\pi_{U_{Q}}\right)=0.$
\item [(iii)] The canonical map $\theta:S\left(V,c\right)\rightarrow\mathrm{gr}^{\prime}\left(U_{Q}\left(V,c,\overline{\beta}\right)\right)$
of Proposition \ref{pro: theta} is an isomorphism of graded braided
bialgebras.
\item [(iv)] $i_{U_{Q}}$
induces an isomorphism between $V$ and $P\left(U_{Q}\left(V,c,\overline{\beta}\right)\right).$ 
\end{enumerate}
\end{thm}
\begin{proof}
Clearly $\emph{B}\left(V,c\right)$ is quadratic if it coincide with
$S_{Q}\left(V,c\right)$. On the other hand, if $\emph{B}\left(V,c\right)$
is quadratic, we get \[
\ker\left(\mathrm{\Gamma}^{T}\right)=\left(\ker\left(\mathrm{\Gamma}_{2}^{T}\right)\right)=\left(\ker\left(c+\mathrm{Id}_{V^{\otimes2}}\right)\right)=\left(E_{2}\left(V,c\right)\right)\]
 and hence $S_{Q}\left(V,c\right)=\emph{B}\left(V,c\right)$. From
this equality, one gets that $P\left(S_{Q}\left(V,c\right)\right)=i_{S}\left(V\right).$
Using this equality the proof is similar to that of \cite[Theorem 6.3]{Ar- Universal}. \end{proof}
\begin{thm}
\label{teo: magnum2}Let $\overline{\beta}:E_{2}\left(V,c\right)\rightarrow V$
be a Qbracket on a braided vector space $\left(V,c\right)$. Assume
that the Nichols algebra $\emph{B}\left(V,c\right)$ is a quadratic
algebra. Then $\left(V,c,\overline{\beta}\right)$ is a QLie algebra
whenever $i_{U_{Q}}:V\rightarrow U_{Q}\left(V,c,\overline{\beta}\right)$
is injective. The converse is true if $\emph{B}\left(V,c\right)$
is Koszul. \end{thm}
\begin{proof}
Let $P:=\left(\mathrm{Id}-\overline{\beta}\right)\left[E_{2}\left(V,c\right)\right].$

Assume that $i_{U_{Q}}:V\rightarrow U_{Q}\left(V,c,\overline{\beta}\right)$
is injective. We have to prove that \eqref{form:Jacobi} holds. Let
$z\in\overline{E_{2}\left(V,c\right)},$ then\begin{eqnarray*}
V & \ni & \overline{\beta}\left(\overline{\beta}_{1}-\overline{\beta}_{2}\right)\left(z\right)=\left(\overline{\beta}-\mathrm{Id}\right)\left(\overline{\beta}_{1}-\overline{\beta}_{2}\right)\left(z\right)+\left(\overline{\beta}_{1}-\mathrm{Id}\right)\left(z\right)+\left(\mathrm{Id}-\overline{\beta}_{2}\right)\left(z\right)\\
 & \in & P+P\otimes V+V\otimes P\subseteq\ker\left(\pi_{U_{Q}}\right).\end{eqnarray*}
By Theorem \ref{teo: magnum}, we have $V\cap\ker\left(\pi_{U_{Q}}\right)=0$
so that $\overline{\beta}\left(\overline{\beta}_{1}-\overline{\beta}_{2}\right)\left(z\right)=0.$

Conversely, assume that $\left(V,c,\overline{\beta}\right)$ is a
QLie algebra and that $\emph{B}\left(V,c\right)$ is Koszul. By Theorem \ref{teo: magnum},
$i_{U_{Q}}:V\rightarrow U_{Q}\left(V,c,\overline{\beta}\right)$ is
injective if and only if the canonical map $\theta:S\left(V,c\right)\rightarrow\mathrm{gr}^{\prime}\left(U_{Q}\left(V,c,\overline{\beta}\right)\right)$
of Proposition \ref{pro: theta} is an isomorphism of graded braided
bialgebras. This means that $P\leq K\oplus V\oplus V^{\otimes2}$
is of PBW type in the sense of \cite{BG}. By \cite[Theorem 0.5]{BG},
$P$ is of PBW if and only if\begin{gather}
P\cap T^{\leq1}=0;\label{form: BG (I)}\\
\left(T^{\leq1}\cdot P\cdot T^{\leq1}\right)\cap T^{\leq2}=P.\label{form: BG (J)}\end{gather}

Let us prove that \eqref{form: BG (I)} is always true. Let $\gamma:=\mathrm{Id}_{V^{\otimes2}}-\overline{\beta}:E_{2}\left(V,c\right)\rightarrow V\oplus V^{\otimes2}.$
Then $\mathrm{Im}\left(\gamma\right)=\left(1-\overline{\beta}\right)\left(E_{2}\left(V,c\right)\right)=P.$
Let $y\in P\cap T^{\leq1}.$ Then $y=\gamma\left(x\right)$ for some
$x\in E_{2}\left(V,c\right)$ so that \[
K\oplus V=T^{\leq1}\ni y=\left(\mathrm{Id}-\overline{\beta}\right)\left(x\right)=x-\overline{\beta}\left(x\right)\in V^{\otimes2}\oplus V.\]
Thus $x=0$ and hence $y=\gamma\left(x\right)=0$ so that \eqref{form: BG (I)}
is always true. Let us check that (\ref{form: BG (J)}) is equivalent to \eqref{form:correctness}
and (\ref{form:Jacobi}). We have\[
P=\left(\mathrm{Id}-\overline{\beta}\right)\left[E_{2}\left(V,c\right)\right]=\left\{ x-\overline{\beta}\left(x\right)\shortmid x\in E_{2}\left(V,c\right)\right\} .\]
The conclusion follows by applying \cite[Lemma 3.3]{BG} to the case
$\alpha=b$ and $\beta=0.$ \end{proof}
\begin{thm}
\label{thm:QLie}Let $\left(V,c,b\right)$ be a braided Lie algebra
such that $\emph{B}\left(V,c\right)$ is a quadratic algebra. Then
$\left(V,c,b^{2}\right)$ is a QLie algebra and $
U\left(V,c,b\right)=U_{Q}\left(V,c,b^{2}\right)$.
\end{thm}
\begin{proof}
By definition \[
U\left(V,c,b\right):=\frac{T\left(V,c\right)}{\left(\left(\mathrm{Id}-b\right)\left[E\left(V,c\right)\right]\right)}\text{\quad and\quad}U_{Q}\left(V,c,b^{2}\right):=\frac{T\left(V,c\right)}{\left(\left(\mathrm{Id}-b^{2}\right)\left[E_{2}\left(V,c\right)\right]\right)}.\]
Since $\left(\mathrm{Id}-b^{2}\right)\left[E_{2}\left(V,c\right)\right]\subseteq\left(\mathrm{Id}-b\right)\left[E\left(V,c\right)\right],$
there is a braided bialgebra homomorphism $\pi:U_{Q}\left(V,c,b^{2}\right)\rightarrow U\left(V,c,b\right)$
which is surjective and such that $\pi\circ i_{U_{Q}}=i_{U}.$ By
Definition \ref{def: c-Lie alg}, we have that $i_{U}$ is injective
so that $i_{U_{Q}}$ is injective too. Clearly $b^{2}:E_{2}\left(V,c\right)\rightarrow V$
is a Qbracket. Hence, by Theorem \ref{teo: magnum}, $i_{U_{Q}}:V\rightarrow U_{Q}\left(V,c,b\right)$
induces an isomorphism between $V$ and $P\left(U_{Q}\left(V,c,b\right)\right).$
Thus, by \cite[Lemma 5.3.3]{Montgomery}, $\pi$ is injective being
injective its restriction to $P\left(U_{Q}\left(V,c,b\right)\right)$.
In other words, $U\left(V,c,b\right)=U_{Q}\left(V,c,b^{2}\right)$.
By Theorem \ref{teo: magnum2}, $\left(V,c,b^{2}\right)$ is a QLie
algebra. 
\end{proof}

\subsection{Lifted quadratic Lie algebras}

We will see that the notion of quadratic Lie algebra naturally leads
to the equivalent notion of lifted quadratic Lie algebra.
\begin{defn}
\label{def:liftedQLie}A \textbf{lifted quadratic Lie algebra} (lifted
QLie algebra for short) is a tern $\left(V,c,\beta\right)$ where
\begin{itemize}
\item $\left(V,c\right)$ is a braided vector space;
\item $\beta:V\otimes V\rightarrow V$ is a map such that\begin{gather}
c\beta_{1}=\beta_{2}c_{1}c_{2},\qquad c\beta_{2}=\beta_{1}c_{2}c_{1},\label{form:bracketbeta}\\
\beta\left(\beta_{1}-\beta_{2}\right)\left(\overline{E_{2}\left(V,c\right)}\right)=0,\label{form:Jacobibeta}\\
\beta c=-\beta,\label{form:antisymmetrybeta}\end{gather}
where \[
\overline{E_{2}\left(V,c\right)}=\left(E_{2}\left(V,c\right)\otimes V\right)\cap\left(V\otimes E_{2}\left(V,c\right)\right)=\left\{ z\in V^{\otimes3}\mid c_{1}\left(z\right)=-z=c_{2}\left(z\right)\right\} .\]
 
\end{itemize}
\end{defn}
\begin{lem}
\label{lem: utile}Let $W$ be a vector space and let $c\in End\left(W\right).$
Assume there are $\alpha,\beta\in K\left[X\right]$ be such that $\alpha\left(c\right)\circ\beta\left(c\right)=0$
and $\gcd\left(\alpha,\beta\right)=1.$ Then
\begin{enumerate}
\item [i)] $\mathrm{Im}\left(\alpha\left(c\right)\right)=\ker\left(\beta\left(c\right)\right);$
\item [ii)] $\mathrm{Im}\left(\alpha\left(c\right)\right)\oplus\mathrm{Im}\left(\beta\left(c\right)\right)=W.$ 
\end{enumerate}
\end{lem}
\begin{proof}
Since GCD$\left(\alpha,\beta\right)=1,$ there are $\alpha^{\prime},\beta^{\prime}\in K\left[X\right]$
such that $\alpha^{\prime}\alpha+\beta^{\prime}\beta=1.$

i) Clearly $\mathrm{Im}\left(\alpha\left(c\right)\right)\subseteq\ker\left(\beta\left(c\right)\right).$
Let $z\in\ker\left(\beta\left(c\right)\right).$ Then \[
z=\left(\alpha^{\prime}\left(c\right)\circ\alpha\left(c\right)+\beta^{\prime}\left(c\right)\circ\beta\left(c\right)\right)\left(z\right)=\alpha\left(c\right)\left[\alpha^{\prime}\left(c\right)\left(z\right)\right]\in\mathrm{Im}\left(\alpha\left(c\right)\right).\]

ii) From $\alpha\alpha^{\prime}+\beta\beta^{\prime}=1$ we deduce
$\mathrm{Im}\left(\alpha\left(c\right)\right)+\mathrm{Im}\left(\beta\left(c\right)\right)=W.$
Let $w\in\mathrm{Im}\left(\alpha\left(c\right)\right)\cap\mathrm{Im}\left(\beta\left(c\right)\right).$
Then $w=\alpha\left(c\right)\left(z\right)$ for some $z\in W,$ whence
\[
z=\alpha^{\prime}\left(c\right)\alpha\left(c\right)\left(z\right)+\beta\left(c\right)\beta^{\prime}\left(c\right)\left(z\right)=\alpha^{\prime}\left(c\right)\left(w\right)+\beta\left(c\right)\beta^{\prime}\left(c\right)\left(z\right)\in\mathrm{Im}\left(\beta\left(c\right)\right).\]
Thus $w=\alpha\left(c\right)\left(z\right)=0$ as $\alpha\left(c\right)\circ\beta\left(c\right)=0.$ 
\end{proof}
\begin{notation} Let $\left(V,c\right)$ be a braided vector space
and assume that $c$ is root of a polynomial $f\in K\left[X\right]$
having $-1$ as a simple root. In this case we will write\[
f=\left(X+1\right)h\]
where $h\in K\left[X\right]$. Note that $\gcd\left(X+1,h\right)=1\ $as
$h\left(-1\right)\neq0$. \end{notation}
\begin{lem}
\label{lem: lifting beta}Let $\left(V,c\right)$ be a braided vector
space. Assume that $c$ is root of a polynomial $f\in K\left[X\right]$
having $-1$ as a simple root.

Then, the assignment $\overline{\beta}\mapsto\overline{\beta}h\left(c\right)$
yields a bijection between the following sets.
\begin{itemize}
\item $\overline{B}=\left\{ \overline{\beta}:E_{2}\left(V,c\right)\rightarrow V\mid\left(V,c,\overline{\beta}\right)\text{ is a QLie algebra}\right\} $.
\item $B=\left\{ \beta:V\otimes V\rightarrow V\mid\left(V,c,\beta\right)\text{ is a lifted QLie algebra}\right\} .$ 
\end{itemize}
\end{lem}
\begin{proof}
Let $\overline{\beta}\in\overline{B}.$ Since, by Lemma \ref{lem: utile},
$E_{2}\left(V,c\right)=\mathrm{Im}\left(h\left(c\right)\right),$
it makes sense to define $\beta:V\otimes V\rightarrow V,\beta\left(z\right)=\overline{\beta}h\left(c\right)\left(z\right).$
Let us check that $\beta\in B.$ For every $z\in V^{\otimes3},$ we
have\[
c\beta_{1}\left(z\right)=c\overline{\beta}_{1}h\left(c\right)_{1}\left(z\right)\overset{\text{(\ref{form:bracket})}}{=}\overline{\beta}_{2}c_{E_{2}\left(V,c\right),V}h\left(c\right)_{1}\left(z\right)=\overline{\beta}_{2}c_{1}c_{2}h\left(c\right)_{1}\left(z\right)\overset{\text{(\ref{ec: braided equation})}}{=}\beta_{2}c_{1}c_{2}\left(z\right).\]
In a similar way one proves that $c\beta_{2}=\beta_{1}c_{2}c_{1}$
holds. Let $z\in\overline{E_{2}\left(V,c\right)}.$ We have\begin{eqnarray*}
\beta\left(\beta_{1}-\beta_{2}\right)\left(z\right) & = & \overline{\beta}h\left(c\right)\left(\overline{\beta}_{1}h\left(c\right)_{1}-\overline{\beta}_{2}h\left(c\right)_{2}\right)\left(z\right)\\
 & = & h\left(-1\right)\overline{\beta}h\left(c\right)\left(\overline{\beta}_{1}-\overline{\beta}_{2}\right)\left(z\right)\\
 & \overset{(\ref{form:correctness})}{=} & h\left(-1\right)^{2}\overline{\beta}\left(\overline{\beta}_{1}-\overline{\beta}_{2}\right)\left(z\right)\\
 & \overset{(\ref{form:Jacobi})}{=} & 0.\end{eqnarray*}
Let $z\in V^{\otimes2},$ we have\[
\beta\left(c+\mathrm{Id}_{V^{\otimes2}}\right)\left(z\right)=\overline{\beta}h\left(c\right)\left(c+\mathrm{Id}_{V^{\otimes2}}\right)\left(z\right)=\overline{\beta}f\left(c\right)\left(z\right)=0.\]
Hence it makes sense to define $\varphi:\overline{B}\rightarrow B,\varphi\left(\overline{\beta}\right)=\beta.$

Conversely, let $\beta\in B.$ Since, by Lemma \ref{lem: utile},
$\mathrm{Im}\left(c+1\right)=\ker\left(h\left(c\right)\right),$ from
(\ref{form:antisymmetrybeta}) we get that there is a unique map $\overline{\beta}:\mathrm{Im}\left(h\left(c\right)\right)=E_{2}\left(V,c\right)\rightarrow V$
such that $\overline{\beta}h\left(c\right)=\beta.$

Since $\gcd\left(h,X+1\right)=1,$ there are $r,s\in K\left[X\right]$
such that $1=hr+s\left(X+1\right).$ Let $z\in\overline{E_{2}\left(V,c\right)}\otimes V=\ker\left(c_{1}+\mathrm{Id}_{V^{\otimes3}}\right).$
We have\begin{eqnarray*}
c\overline{\beta}_{1}\left(z\right) & = & c\overline{\beta}_{1}\left[h\left(c\right)_{1}r\left(c\right)_{1}+s\left(c\right)_{1}\left(c_{1}+\mathrm{Id}_{V^{\otimes3}}\right)\right]\left(z\right)\\
 & = & c\overline{\beta}_{1}h\left(c\right)_{1}r\left(c\right)_{1}\left(z\right)=c\beta_{1}r\left(c\right)_{1}\left(z\right)\overset{\text{\text{(\ref{form:bracketbeta})}}}{=}\beta_{2}c_{1}c_{2}r\left(c\right)_{1}\left(z\right)\\
 & = & \overline{\beta}_{2}h\left(c\right)_{2}c_{1}c_{2}r\left(c\right)_{1}\left(z\right)\overset{\text{(\ref{ec: braided equation})}}{=}\overline{\beta}_{2}c_{1}c_{2}h\left(c\right)_{1}r\left(c\right)_{1}\left(z\right)\\
 & = & \overline{\beta}_{2}c_{1}c_{2}\left[h\left(c\right)_{1}r\left(c\right)_{1}+s\left(c\right)_{1}\left(c_{1}+\mathrm{Id}_{V^{\otimes3}}\right)\right]\left(z\right)=\overline{\beta}_{2}c_{1}c_{2}\left(z\right)\\
 &=&\overline{\beta}_{2}c_{E_{2}\left(V,c\right),V}\left(z\right)\end{eqnarray*}
Similarly one proves that $c\overline{\beta}_{2}=\overline{\beta}_{1}c_{V,E_{2}\left(V,c\right)}.$
Let $z\in\overline{E_{2}\left(V,c\right)}.$ We have\begin{eqnarray*}
0 & \overset{(\ref{form:Jacobibeta})}{=} & \beta\left(\beta_{1}-\beta_{2}\right)\left(z\right)\\
 & = & \overline{\beta}h\left(c\right)\left(\overline{\beta}_{1}h\left(c\right)_{1}-\overline{\beta}_{2}h\left(c\right)_{2}\right)\left(z\right)\\
 & = & h\left(-1\right)\overline{\beta}h\left(c\right)\left(\overline{\beta}_{1}-\overline{\beta}_{2}\right)\left(z\right)\\
 & \overset{(\ref{form:correctness})}{=} & h\left(-1\right)^{2}\overline{\beta}\left(\overline{\beta}_{1}-\overline{\beta}_{2}\right)\left(z\right).\end{eqnarray*}
Since $h\left(-1\right)\neq0,$ we get $\overline{\beta}\left(\overline{\beta}_{1}-\overline{\beta}_{2}\right)\left(z\right)=0.$ Hence it makes sense to define $\psi:B\rightarrow\overline{B},\psi\left(\beta\right)=\overline{\beta}.$ The uniqueness of $\overline{\beta}$ implies that $\varphi$ and
$\psi$ are mutual inverses. \end{proof}
\begin{rem}
\label{rem:UQ}With notations of Lemma \ref{lem: lifting beta}, if
$\overline{\beta}\in\overline{B},$ we have \begin{eqnarray*}
U_{Q}\left(V,c,\overline{\beta}\right) & = & \frac{T\left(V,c\right)}{\left(\left(\mathrm{Id}-\overline{\beta}\right)\left[E_{2}\left(V,c\right)\right]\right)}\\
 & = & \frac{T\left(V,c\right)}{\left(\left(\mathrm{Id}-\overline{\beta}\right)h\left(c\right)\left(z\right)\mid z\in V^{\otimes2}\right)}\\
 & = & \frac{T\left(V,c\right)}{\left(h\left(c\right)\left(z\right)-\beta\left(z\right)\mid z\in V^{\otimes2}\right)}.\end{eqnarray*}
 \end{rem}
\begin{defn}
The previous remark justifies the following notation:\[
U_{Q}\left(V,c,\beta\right):=\frac{T\left(V,c\right)}{\left(h\left(c\right)\left(z\right)-\beta\left(z\right)\mid z\in V^{\otimes2}\right)}.\]
\end{defn}
\begin{thm}
\label{thm:quadratic}Let $A$ be a primitively generated connected
braided bialgebra and let $(P,c)$ be the braided vector space of
primitive elements of $A$. Assume that the Nichols algebra $\emph{B}\left(P,c\right)$
is a quadratic algebra. We have that
\begin{enumerate}
\item [1)] if $\overline{\beta}$ is the restriction of $m_A$ to $E_2(P,c)$, then $\left(P,c,\overline{\beta}\right)$ is a QLie algebra and $A$ is isomorphic to the universal enveloping algebra $U_{Q}(P,c,\overline{\beta})$;
\item [2)] if $c$ is root of a polynomial $f\in K\left[X\right]$ of the form $f=(X+1)h$ where $h(-1)\neq 0$ and $\beta=m_Ah(c)$, then $\left(P,c,\beta\right)$ is a lifted
QLie algebra and $A$ is isomorphic to the
universal enveloping algebra $U_{Q}\left(V,c,\beta\right)$.
\end{enumerate}
\end{thm}
\begin{proof}
1) Let $(P,c,b)$ be the infinitesimal braided Lie algebra of $A$.
Then, by \cite[Theorems 9.3 and 6.5]{Ar- Universal}, we have that
$A$ is isomorphic to the universal enveloping algebra $U(P,c,b)$
as a braided bialgebra, where $b(z)=m_A^{t-1}(z)$ for all $z\in E_t(P,c)$. Set $\overline{\beta}:=b^{2}$. Then $\overline{\beta}$ is the restriction of $m_A$ to $E_2(P,c)$. By Theorem
\ref{thm:QLie}, $\left(P,c,\overline{\beta}\right)$ is a QLie algebra and $U\left(P,c,b\right)=U_{Q}\left(P,c,\overline{\beta}\right).$

2) By 1) and Lemma \ref{lem: lifting beta}, we can consider the lifted
QLie algebra $\left(V,c,\beta\right)$ where $\beta=\overline{\beta}h(c)$. By Remark \ref{rem:UQ}, $U_{Q}\left(P,c,\overline{\beta}\right)=U_{Q}\left(V,c,\beta\right)$.\end{proof}
\begin{lem}
Let $\left(V,c\right)$ be a braided vector space and let $\beta:V\otimes V\rightarrow V$
satisfying (\ref{form:bracketbeta}) and (\ref{form:antisymmetrybeta}).
Then\[
\beta\left(\beta_{1}+\beta_{2}\right)\left(\overline{E_{2}\left(V,c\right)}\right)=0.\]
Therefore, if $char\left(K\right)=2$ then (\ref{form:Jacobibeta})
holds whence $\left(V,c,\beta\right)$ is a lifted QLie algebra.
If $char\left(K\right)\neq2,$ then (\ref{form:Jacobibeta}) holds
if and only if $\beta\beta_{1}\left(\overline{E_{2}\left(V,c\right)}\right)=0$
if and only if $\beta\beta_{2}\left(\overline{E_{2}\left(V,c\right)}\right)=0$. \end{lem}
\begin{proof}
Let $z\in\overline{E_{2}\left(V,c\right)}.$ Then\[
\beta\beta_{1}\left(z\right)=\beta\beta_{1}c_{2}c_{1}\left(z\right)\overset{(\ref{form:bracketbeta})}{=}\beta c\beta_{2}\left(z\right)\overset{(\ref{form:antisymmetrybeta})}{=}-\beta\beta_{2}\left(z\right).\]
 
\end{proof}

\subsection{Categorical subspaces}
\begin{defn}
\cite[2.2]{Kharchenko- connected} A subspace $L$ of a braided vector
space $\left(V,c\right)$ is said to be \textbf{categorical} if \begin{equation}
c\left(L\otimes V\right)\subseteq V\otimes L\qquad\text{and}\qquad c\left(V\otimes L\right)\subseteq L\otimes V.\label{form: categorical}\end{equation}
 \end{defn}
\begin{lem}
\label{lem: categorical}Let $\left(V,c\right)$ be a braided vector
space and let $\beta:V\otimes V\rightarrow V$ be such that (\ref{form:bracketbeta})
holds. Then $L=\mathrm{Im}\left(\beta\right)$ is a categorical subspace
of $V$. Furthermore, if $\left(V,c,\beta\right)$ is a lifted QLie
algebra, then $\left(L,c_{L},\beta_{L}\right)$ is a lifted QLie subalgebra
of $\left(V,c,\beta\right)$ where $c_{L}$ and $\beta_{L}$ are the
obvious restrictions of $c$ and $\beta$ respectively. \end{lem}
\begin{proof}
We have that\begin{eqnarray*}
c\beta_{1} & = & \beta_{2}c_{1}c_{2}\Rightarrow c\left(L\otimes V\right)\subseteq V\otimes L,\\
c\beta_{2} & = & \beta_{1}c_{2}c_{1}\Rightarrow c\left(V\otimes L\right)\subseteq L\otimes V.\end{eqnarray*}
Assume that $\left(V,c,\beta\right)$ is a lifted QLie algebra. Then,
by the foregoing $c\left(L\otimes L\right)\subseteq L\otimes L$ so
that $c_{L}:L\otimes L\rightarrow L\otimes L$ is well defined. The
conclusion follows. \end{proof}
\begin{prop}
\label{pro: f(-1)}Let $\left(V,c,\overline{\beta}\right)$ be a QLie
algebra with $V$ finite dimensional. Then $c$ has a minimal polynomial
$f$. Moreover if $f\left(-1\right)\neq0$ then $U_{Q}\left(V,c,\overline{\beta}\right)=T\left(V,c\right).$ \end{prop}
\begin{proof}
The first part is clear. Assume $f\left(-1\right)\neq0.$ Then $\gcd\left(f,X+1\right)=1$
so that $c+\mathrm{Id}_{V^{\otimes2}}$ is invertible by Bézout identity.
Therefore $E_{2}\left(V,c\right)=\ker\left(c+\mathrm{Id}_{V^{\otimes2}}\right)=0$
whence $U_{Q}\left(V,c,\overline{\beta}\right)=T\left(V,c\right)$.
\end{proof}

Next proposition will be used to normalize $\beta$ if necessary.
\begin{prop}
\label{pro: normalization of beta}Let $\left(V,c,\beta\right)$ be
a lifted QLie algebra and let $\lambda\in K.$ Then $\left(V,c,\lambda\beta\right)$
is a lifted QLie algebra. Moreover $\lambda\mathrm{Id}_{V}:\left(V,c,\lambda\beta\right)\rightarrow\left(V,c,\beta\right)$
is a morphism of lifted QLie algebras which is an isomorphism whenever
$\lambda\neq0.$ \end{prop}
\begin{proof}
It is straightforward. 
\end{proof}

\section{Low dimensional cases\label{sec:Low-dimensional-cases}}

The main aim of this section is to classify lifted QLie algebras $\left(V,c,\beta\right)$
such that $\mathrm{dim}\left(V\right)\leq2$ and $\beta\neq0$, under
mild assumptions. We will complete the classification in case $\mathrm{dim}\left(\mathrm{Im}\beta\right)=1$.
The case $\mathrm{dim}\left(\mathrm{Im}\beta\right)=2$ will be partially
treated in the Appendix \ref{sec:appendix}.

\subsection{Dimension 1}

In this subsection we treat the case $\mathrm{dim}\left(V\right)=1$.
\begin{prop}
\label{pro: Lie of dim 1}Assume $char\left(K\right)\neq2$ and let
$\left(V,c,\beta\right)$ be a one dimensional lifted QLie algebra.
Then $\beta=0.$ \end{prop}
\begin{proof}
Let $V=Kx.$ Then $\beta\left(x\otimes x\right)=\lambda x$ for some
$\lambda\in K.$ Assume $\lambda\neq0.$ Now $c\left(x\otimes x\right)=\gamma x\otimes x$
for some $\gamma\in K.$ We have \[
-\lambda x=-\beta\left(x\otimes x\right)\overset{(\ref{form:antisymmetrybeta})}{=}\beta c\left(x\otimes x\right)=\gamma\beta\left(x\otimes x\right)=\gamma\lambda x.\]
From this we get $\lambda\left(\gamma+1\right)=0$ which implies $\gamma=-1$
whence $c=-\mathrm{Id}_{V^{\otimes2}}.$ Then, from (\ref{form:bracketbeta}),
we have $-\beta_{1}=\beta_{2}$ which implies \[
-\lambda x\otimes x=-\beta_{1}\left(x\otimes x\otimes x\right)=\beta_{2}\left(x\otimes x\otimes x\right)=\lambda x\otimes x.\]
Since $char\left(K\right)\neq2$, we get $\lambda=0,$ a contradiction. 
\end{proof}

\subsection{Dimension 2}

In this subsection we treat the case $\mathrm{dim}\left(V\right)=2$.
\begin{cor}
\label{coro: b11=00003D0}Assume $char\left(K\right)\neq2$ and let
$\left(V,c,\beta\right)$ be a lifted QLie algebra. If $L:=\mathrm{Im}\beta$
is one dimensional, then $\beta\left(L\otimes L\right)=0.$ \end{cor}
\begin{proof}
By Lemma \ref{lem: categorical}, $\left(L,c_{L},\beta_{L}\right)$
is a lifted QLie subalgebra of $\left(V,c,\beta\right).$ By Proposition
\ref{pro: Lie of dim 1}, we have $\beta_{L}=0.$ 
\end{proof}

Now we introduce some notation that will be needed afterwards.

\begin{notation} \label{not: dim 2} We take the following assumtions and notations.
\begin{itemize}
\item $K$ is a field such that $char\left(K\right)\neq2.$
\item $\left(V,c\right)$ is a two dimensional braided vector space. Chosen a
basis $x_{1},x_{2}$ for $V$ and taken $x_{1}\otimes x_{1},x_{2}\otimes x_{1},x_{1}\otimes x_{2},x_{2}\otimes x_{2}$
as a basis fo $V\otimes V$, we will write\[
c=\begin{pmatrix}c_{11}^{11} & c_{11}^{21} & c_{11}^{12} & c_{11}^{22}\\
c_{21}^{11} & c_{21}^{21} & c_{21}^{12} & c_{21}^{22}\\
c_{12}^{11} & c_{12}^{21} & c_{12}^{12} & c_{12}^{22}\\
c_{22}^{11} & c_{22}^{21} & c_{22}^{12} & c_{22}^{22}\end{pmatrix}.\]
 i.e. $c\left(x_{i}\otimes x_{j}\right)=c_{11}^{ij}x_{1}\otimes x_{1}+c_{21}^{ij}x_{2}\otimes x_{1}+c_{12}^{ij}x_{1}\otimes x_{2}+c_{22}^{ij}x_{2}\otimes x_{2}$.
\item $\beta:V\otimes V\rightarrow V$ is a nonzero map such that $\left(V,c,\beta\right)$
is a lifted QLie algebra. With the same basis as above we will write\[
\beta=\begin{pmatrix}\beta_{1}^{11} & \beta_{1}^{21} & \beta_{1}^{12} & \beta_{1}^{22}\\
\beta_{2}^{11} & \beta_{2}^{21} & \beta_{2}^{12} & \beta_{2}^{22}\end{pmatrix}\]
 i.e. $\beta\left(x_{i}\otimes x_{j}\right)=\beta_{1}^{ij}x_{1}+\beta_{2}^{ij}x_{2}$. 
 \item  $c$ has minimal polynomial $f\in K\left[X\right]$ of the form $f=\left(X+1\right)h$ for some $h\in K\left[X\right]$ such that
$h\left(-1\right)\neq0.$ 
\end{itemize}
\end{notation}

Inside the proof of next result, several computations have been handled with the help of the Computer Algebra System AXIOM \cite{Axiom}. 

\begin{thm}
\label{thm:teoClassif}Assume $\dim\left(\mathrm{Im}\beta\right)=1.$
Then, there exist a basis $x_{1},x_{2}$ for $V$ and $\gamma\in K$
such that $c$, $\beta,$ $U_{Q}\left(V,c,\beta\right)$ and $f$
take only one of the canonical forms (up to an isomorphism of lifted
QLie algebras) in Table \ref{tab: Im(Beta)=00003D1}. Moreover, if $c$ and $\beta$ are as in any of the cases of Table
\ref{tab: Im(Beta)=00003D1}, then $\left(V,c,\beta\right)$ is a
lifted QLie algebra.
\begin{table}[tbh]
  \setlength\arraycolsep{3pt}
  \setlength\tabcolsep{2pt}
{\small \[
\begin{array}{|c|c|c|c|c|}
\hline \strut & \strut & \strut & \strut & \strut\\
f & c & \beta & U_{Q}\left(V,c,\beta\right) & \gamma\\
\strut & \strut & \strut & \strut & \strut\\
\hline \strut & \strut & \strut & \strut & \strut\\
X^{2}-1 & 
\begin{pmatrix} 1 & 0 & 0 & 0\\
0 & 0 & 1 & 0\\
0 & 1 & 0 & 0\\
0 & 0 & 0 & 1\end{pmatrix} & \begin{pmatrix}0 & 1 & -1 & 0\\
0 & 0 & 0 & 0\end{pmatrix} & \frac{\strut T\left(V,c\right)}{\strut(x_{2}x_{1}-x_{1}x_{2}+x_{1})} & \strut\\
\strut & \strut & \strut & \strut & \strut\\
X^{2}-1 & \begin{pmatrix}1 & 1 & -1 & 0\\
0 & 0 & 1 & 0\\
0 & 1 & 0 & 0\\
0 & 0 & 0 & 1\end{pmatrix} & \begin{pmatrix}0 & 1 & -1 & 0\\
0 & 0 & 0 & 0\end{pmatrix} & \frac{\strut T\left(V,c\right)}{\left(\left(x_{1}\right)^{2}-x_{2}x_{1}+x_{1}x_{2}-x_{1}\right)} & \strut\\
\strut & \strut & \strut & \strut & \strut\\
X^{2}-1 & \begin{pmatrix}-1 & 0 & 0 & \gamma\\
0 & 0 & 1 & 0\\
0 & 1 & 0 & 0\\
0 & 0 & 0 & 1\end{pmatrix} & \begin{pmatrix}0 & 1 & -1 & 0\\
0 & 0 & 0 & 0\end{pmatrix} & \frac{\strut T\left(V,c\right)}{\left(\begin{array}{c}
{\scriptstyle (x_{1})^{2},}\\
{\scriptstyle x_{2}x_{1}-x_{1}x_{2}+x_{1}}\end{array}\right)} & \begin{array}{c}
0,1\text{ or }\\
\sqrt{\gamma}\notin K\end{array}\\
\strut & \strut & \strut & \strut & \strut\\
X^{2}-1 & \begin{pmatrix}1 & 0 & 0 & \gamma\\
0 & 0 & 1 & 0\\
0 & 1 & 0 & 0\\
0 & 0 & 0 & -1\end{pmatrix} & \begin{pmatrix}0 & 0 & 0 & 1\\
0 & 0 & 0 & 0\end{pmatrix} & \frac{\strut T\left(V,c\right)}{\left(\begin{array}{c}
{\scriptstyle x_{2}x_{1}-x_{1}x_{2},}\\
{\scriptstyle \gamma\left(x_{1}\right)^{2}-2\left(x_{2}\right)^{2}-x_{1}}\end{array}\right)} & \begin{array}{c}
0,1\text{ or }\\
\sqrt{\gamma}\notin K\end{array}\\
\strut & \strut & \strut & \strut & \strut\\
(X^{2}-1)X & \begin{pmatrix}0 & 1 & 0 & 0\\
0 & 0 & 1 & 0\\
0 & 1 & 0 & 0\\
0 & 0 & 0 & 1\end{pmatrix} & \begin{pmatrix}0 & 1 & -1 & 0\\
0 & 0 & 0 & 0\end{pmatrix} & \frac{\strut T\left(V,c\right)}{\left(\left(x_{1}\right)^{2}-x_{2}x_{1}+x_{1}x_{2}+x_{1}\right)} & \strut\\
\strut & \strut & \strut & \strut & \strut\\
(X^{2}-1)X & \begin{pmatrix}0 & 0 & 0 & 0\\
0 & 0 & \gamma & 0\\
0 & \gamma^{-1} & 0 & 0\\
0 & 0 & 0 & 1\end{pmatrix} & \begin{pmatrix}0 & 1 & -\gamma & 0\\
0 & 0 & 0 & 0\end{pmatrix} & \frac{\strut T\left(V,c\right)}{\strut(-\gamma x_{2}x_{1}+x_{1}x_{2}+\gamma x_{1})} & \gamma\neq0,1\\
\strut & \strut & \strut & \strut & \strut\\
(X^{2}-1)(X-\gamma) & \begin{pmatrix}\gamma & 0 & 0 & 0\\
0 & 0 & 1 & 0\\
0 & 1 & 0 & 0\\
0 & 0 & 0 & 1\end{pmatrix} & \begin{pmatrix}0 & 1 & -1 & 0\\
0 & 0 & 0 & 0\end{pmatrix} & \frac{\strut T\left(V,c\right)}{\big(\left(1+\gamma\right)\left(x_{2}x_{1}-x_{1}x_{2}\right)-x_{1}\big)} & \gamma\neq\pm1\\
\strut & \strut & \strut & \strut & \strut\\
(X^{2}-1)(X-1) & \begin{pmatrix}1 & 0 & 0 & \gamma\\
0 & 0 & 1 & 0\\
0 & 1 & 0 & 0\\
0 & 0 & 0 & 1\end{pmatrix} & \begin{pmatrix}0 & 1 & -1 & 0\\
0 & 0 & 0 & 0\end{pmatrix} & \frac{\strut T\left(V,c\right)}{\big(2\left(x_{2}x_{1}-x_{1}x_{2}\right)-x_{1}\big)} & \begin{array}{c}
1\text{ or }\\
\sqrt{\gamma}\notin K\end{array}\\
\strut & \strut & \strut & \strut & \strut\\\hline \end{array}\]
} 

\caption{Lifted two dimensional QLie algebras with $\dim\left(\mathrm{Im}\beta\right)=1.$}

\label{tab: Im(Beta)=00003D1} 
\end{table}
 \end{thm}
\begin{rem}
The first case is the unique ordinary two-dimensional Lie algebra
in $char\left(K\right)\neq2$. \end{rem}
\begin{proof}
First we will show that we can reduce to one of the cases in Table
\ref{tab: Im(Beta)=00003D1}.

By hypothesis $\dim\left(\mathrm{Im}\beta\right)=1.$ Choose $x_{1}\in V$
such that $L:=\mathrm{Im}\beta=Kx_{1}.$ Complete $x_{1}$ to a basis
$x_{1},x_{2}$ of $V.$ By Lemma \ref{lem: categorical}, $L$ is
a categorical subspace of $V$ so that $c\left(L\otimes V\right)\subseteq V\otimes L,c\left(V\otimes L\right)\subseteq L\otimes V$
and $c\left(L\otimes L\right)\subseteq L\otimes L.$ Moreover, by
Corollary \ref{coro: b11=00003D0}, we have $\beta_{1}^{11}=0.$ These
facts imply that $c$ and $\beta$ take the following form:\begin{equation}
c=\begin{pmatrix}c_{11}^{11} & c_{11}^{21} & c_{11}^{12} & c_{11}^{22}\\
0 & 0 & c_{21}^{12} & c_{21}^{22}\\
0 & c_{12}^{21} & 0 & c_{12}^{22}\\
0 & 0 & 0 & c_{22}^{22}\end{pmatrix},\qquad\beta=\begin{pmatrix}0 & \beta^{21} & \beta^{12} & \beta^{22}\\
0 & 0 & 0 & 0\end{pmatrix}.\label{form:normalized1}\end{equation}
We need the following result.
\begin{lem}
\label{lem: c1111 non zero}We have $-\beta^{21}=\beta^{12}c_{12}^{21}$
and $-\beta^{12}=c_{21}^{12}\beta^{21}.$ Moreover, if $c_{11}^{11}\neq0,$
then $\beta^{21}+\beta^{12}=0.$ \end{lem}
\begin{proof}
We have \begin{eqnarray*}
 &  & -\beta^{21}x_{1}\overset{(\ref{form:normalized1})}{=}-\beta\left(x_{2}\otimes x_{1}\right)\overset{(\ref{form:antisymmetrybeta})}{=}\beta c\left(x_{2}\otimes x_{1}\right)\overset{(\ref{form:normalized1})}{=}c_{12}^{21}\beta^{12}x_{1},\\
 &  & -\beta^{12}x_{1}\overset{(\ref{form:normalized1})}{=}-\beta\left(x_{1}\otimes x_{2}\right)\overset{(\ref{form:antisymmetrybeta})}{=}\beta c\left(x_{1}\otimes x_{2}\right)\overset{(\ref{form:normalized1})}{=}c_{21}^{12}\beta^{21}x_{1}\end{eqnarray*}
so that $-\beta^{21}=\beta^{12}c_{12}^{21}$ and $-\beta^{12}=c_{21}^{12}\beta^{21}.$
Moreover
\begin{equation*}
  \begin{split}
    \beta^{12}c_{11}^{11}c_{12}^{21}x_{1}\otimes x_{1}
     &\overset{(\ref{form:normalized1})}{=}\beta_{2}c_{1}c_{2}\left(x_{1}\otimes x_{2}\otimes x_{1}\right)\\
     &\overset{(\ref{form:bracketbeta})}{=}c\beta_{1}\left(x_{1}\otimes x_{2}\otimes x_{1}\right)\\
     &\overset{(\ref{form:normalized1})}{=}c_{11}^{11}\beta^{12}x_{1}\otimes x_{1}.
  \end{split}
\end{equation*}
If $c_{11}^{11}\neq0,$ we get $\beta^{12}c_{12}^{21}=\beta^{12}.$ 
\end{proof}
CASE 1) Assume $c_{11}^{11}\neq0.$

By Lemma \ref{lem: c1111 non zero}, $c$ and $\beta$ take the following
form:\begin{equation}
c=\begin{pmatrix}c_{11}^{11} & c_{11}^{21} & c_{11}^{12} & c_{11}^{22}\\
0 & 0 & c_{21}^{12} & c_{21}^{22}\\
0 & c_{12}^{21} & 0 & c_{12}^{22}\\
0 & 0 & 0 & c_{22}^{22}\end{pmatrix},\qquad\beta=\begin{pmatrix}0 & \beta^{21} & -\beta^{21} & \beta^{22}\\
0 & 0 & 0 & 0\end{pmatrix}.\label{form:normalized2}\end{equation}
CASE 1.1) Assume $\beta^{21}=0.$ Then, since $\dim\left(\mathrm{Im}\beta\right)=1,$
we have $\beta^{22}\neq0$ so that, by Proposition \ref{pro: normalization of beta},
we can assume $\beta^{22}=1.$

Using (\ref{form:antisymmetrybeta}), we get $c_{22}^{22}=-1.$

Using the left-hand side of (\ref{form:bracketbeta}), we get $c_{21}^{12}=1.$

Using the right-hand side of (\ref{form:bracketbeta}), we get $c_{11}^{11}=1,c_{12}^{21}=1,c_{11}^{21}=0.$

Using the left-hand side of (\ref{form:bracketbeta}), we get $c_{11}^{12}=0.$

Using (\ref{ec: braided equation}), we get $c_{12}^{22}=c_{21}^{22}.$
In conclusion we have\begin{equation}
c=\begin{pmatrix}1 & 0 & 0 & c_{11}^{22}\\
0 & 0 & 1 & c_{21}^{22}\\
0 & 1 & 0 & c_{21}^{22}\\
0 & 0 & 0 & -1\end{pmatrix},\qquad\beta=\begin{pmatrix}0 & 0 & 0 & 1\\
0 & 0 & 0 & 0\end{pmatrix}.\label{form:normalized3}\end{equation}
With the basis change in $V$ given by $x_{1}=tx_{1}^{\prime},x_{2}=c_{21}^{22}x_{1}^{\prime}/2+x_{2}^{\prime}$
we arrive at\[
c=\begin{pmatrix}1 & 0 & 0 & t^{2}g\\
0 & 0 & 1 & 0\\
0 & 1 & 0 & 0\\
0 & 0 & 0 & -1\end{pmatrix},\qquad\beta=\begin{pmatrix}0 & 0 & 0 & 1\\
0 & 0 & 0 & 0\end{pmatrix}.\]
for $g=\left(c_{21}^{22}\right)^{2}/2+c_{11}^{22}.$ Hence we are
in the fourth case in Table \ref{tab: Im(Beta)=00003D1}. In fact,
if $g\neq0,1$ and $\sqrt{g}\in K$ then one can choose $t=1/\sqrt{g}$.
Otherwise choose $t=1$ and $\gamma=g.$

CASE 1.2) Assume $\beta^{21}\neq0.$ Then, by Proposition \ref{pro: normalization of beta},
we can assume $\beta^{21}=1.$

Using (\ref{form:antisymmetrybeta}), we get $c_{12}^{21}=1=c_{21}^{12}.$

Using the left-hand side of (\ref{form:bracketbeta}), we get $c_{12}^{22}=0,c_{22}^{22}=1.$

Using the right-hand side of (\ref{form:bracketbeta}), we get $c_{21}^{22}=0$.

Using (\ref{form:antisymmetrybeta}), we get $\beta^{22}=0.$ In conclusion
we have \begin{equation}
c=\begin{pmatrix}c_{11}^{11} & c_{11}^{21} & c_{11}^{12} & c_{11}^{22}\\
0 & 0 & 1 & 0\\
0 & 1 & 0 & 0\\
0 & 0 & 0 & 1\end{pmatrix},\qquad\beta=\begin{pmatrix}0 & 1 & -1 & 0\\
0 & 0 & 0 & 0\end{pmatrix}.\label{form:normalized4}\end{equation}
We point out that $c$ and $\beta$ as in (\ref{form:normalized4})
always fulfill (\ref{form:antisymmetrybeta}) and (\ref{form:bracketbeta}).

CASE 1.2.1) Assume $c_{11}^{21}\neq c_{11}^{12}.$

Using (\ref{ec: braided equation}), we get $c_{11}^{22}=0,c_{11}^{11}=1,c_{11}^{12}=-c_{11}^{21}$
(this implies $c_{11}^{21}\neq0$ by assumption).

With the basis change in $V$ given by $x_{1}=x_{1}^{\prime}/c_{11}^{21},x_{2}=x_{2}^{\prime}$
we arrive at the second case in Table \ref{tab: Im(Beta)=00003D1}.

CASE 1.2.2) Assume $c_{11}^{21}=c_{11}^{12}.$

CASE 1.2.2.1) Assume $c_{11}^{11}=1.$

Using (\ref{ec: braided equation}), we get $c_{11}^{21}=0.$

CASE 1.2.2.1.1) If $c_{11}^{22}=0,$ we are in the first case in Table
\ref{tab: Im(Beta)=00003D1}.

CASE 1.2.2.1.2) If $c_{11}^{22}\neq0,$ we are in the eighth case
in Table \ref{tab: Im(Beta)=00003D1}. In fact, if $c_{11}^{22}\neq1$
and $\sqrt{c_{11}^{22}}\in K$ then, with the basis change in $V$
given by $x_{1}=x_{1}^{\prime}/\sqrt{c_{11}^{22}},x_{2}=x_{2}^{\prime},$
we can reduce to the case $c_{11}^{22}=1.$

CASE 1.2.2.2) Assume $c_{11}^{11}=-1.$

With the basis change in $V$ given by $x_{1}=x_{1}^{\prime},x_{2}=-x_{1}^{\prime}c_{11}^{21}/2+x_{2}^{\prime}$,
we can reduce to the case $c_{11}^{21}=c_{11}^{12}=0.$ Now we are
in the third case in Table \ref{tab: Im(Beta)=00003D1}. In fact,
if $c_{11}^{22}\neq0,1$ and $\sqrt{c_{11}^{22}}\in K$ then, with
the basis change in $V$ given by $x_{1}=x_{1}^{\prime}/\sqrt{c_{11}^{22}},x_{2}=x_{2}^{\prime},$
we can reduce to the case $c_{11}^{22}=1.$

CASE 1.2.2.3) Assume $c_{11}^{11}\neq\pm1.$ Using (\ref{ec: braided equation}),
we get $c_{11}^{22}=\left(c_{11}^{21}\right)^{2}/\left(c_{11}^{11}-1\right).$
With the basis change in $V$ given by $x_{1}=\left(c_{11}^{11}-1\right)x_{1}^{\prime},x_{2}=x_{1}^{\prime}c_{11}^{21}+x_{2}^{\prime}$
we reduce to the seventh case in Table \ref{tab: Im(Beta)=00003D1}.

CASE 2) Assume $c_{11}^{11}=0.$ By Lemma \ref{lem: c1111 non zero},
we have $\ $\begin{equation}
\beta^{12}=0\Leftrightarrow\beta^{21}=0.\label{form: utile}\end{equation}
%One can check that $c_{22}^{22}\neq 1$ implies $\beta =0$ contradicting the
%assumption $\dim \left( \mathrm{Im}\beta \right) =1.$

CASE 2.1) $c_{22}^{22}\neq1.$ Let us prove that this case does not
occur.

CASE 2.1.1) $\beta^{21}=0.$ By (\ref{form: utile}) we get $\beta^{12}=0.$
In this case, since $\dim\left(\mathrm{Im}\beta\right)=1,$ by Proposition
\ref{pro: normalization of beta}, we can assume $\beta^{22}=1.$

Using (\ref{form:antisymmetrybeta}), we get $c_{22}^{22}=-1.$

Using the left-hand side of (\ref{form:bracketbeta}), we get $c_{12}^{21}=0$
and $c_{21}^{12}=1.$

Using the right-hand side of (\ref{form:bracketbeta}), we get a contradiction.

CASE 2.1.2) $\beta^{21}\neq0.$ By Proposition \ref{pro: normalization of beta},
we can assume $\beta^{21}=1.$

Using the left-hand side of (\ref{form:bracketbeta}), we get $c_{21}^{12}=0.$

Using (\ref{form:antisymmetrybeta}), we get a contradiction.

CASE 2.2) $c_{22}^{22}=1.$

CASE 2.2.1) $\beta^{22}\neq0.$ By Proposition \ref{pro: normalization of beta},
we can assume $\beta^{22}=1.$

Using (\ref{form:antisymmetrybeta}) and (\ref{form: utile}) we obtain
$\beta^{21}\neq0$ and $\beta^{12}\neq0.$

CASE 2.2.1.1) $\beta^{12}=-\beta^{21}.$

Using (\ref{form:antisymmetrybeta}), we get $c_{12}^{21}=1=c_{21}^{12}.$

Using the left-hand side of (\ref{form:bracketbeta}), we get a contradiction.

CASE 2.2.1.2) $\beta^{12}+\beta^{21}\neq0.$

With the basis change in $V$ given by $x_{1}=\beta^{12}\beta^{21}\left(\beta^{12}+\beta^{21}\right)x_{1}^{\prime},x_{2}=\beta^{12}\beta^{21}x_{1}^{\prime}+x_{2}^{\prime}$
we reduce to the following CASE 2.2.2).

CASE 2.2.2) $\beta^{22}=0.$ By Proposition \ref{pro: normalization of beta},
we can assume $\beta^{21}=1.$

CASE 2.2.2.1) $\beta^{12}=-1.$

Using (\ref{form:antisymmetrybeta}), we get $c_{12}^{21}=1=c_{21}^{12}$
and $c_{21}^{22}=c_{12}^{22}.$

Using the left-hand side of (\ref{form:bracketbeta}), we get $c_{12}^{22}=0.$

Using (\ref{ec: braided equation}), we get $c_{11}^{22}=-c_{11}^{12}c_{11}^{21}.$

CASE 2.2.2.1.1) $c_{11}^{12}=c_{11}^{21}.$ With the basis change
in $V$ given by $x_{1}=x_{1}^{\prime},x_{2}=-x_{1}^{\prime}c_{11}^{21}+x_{2}^{\prime}$
we reduce to the seventh case in Table \ref{tab: Im(Beta)=00003D1}
for $\gamma=0$.

CASE 2.2.2.1.2) $c_{11}^{12}\neq c_{11}^{21}.$ With the basis change
in $V$ given by 
\begin{equation*}
  \begin{split}
    x_{1}&=x_{1}^{\prime}/\left(c_{11}^{21}-c_{11}^{12}\right)\\
    x_{2}&=-x_{1}^{\prime}c_{11}^{12}/\left(c_{11}^{21}-c_{11}^{12}\right)+x_{2}^{\prime}
  \end{split}
\end{equation*}
we are in the fifth case in Table \ref{tab: Im(Beta)=00003D1}.

CASE 2.2.2.2) $\beta^{12}\neq-1.$

Using (\ref{form:antisymmetrybeta}), we get $c_{12}^{21}\neq0.$

Using the left-hand side of (\ref{form:bracketbeta}), we get $c_{11}^{21}=0=c_{21}^{22}.$

Using (\ref{form:antisymmetrybeta}) and (\ref{form: utile}), we
get $c_{12}^{22}=0$ and hence $c_{21}^{12}=-\beta^{12}\neq0.$

Using the right-hand side of (\ref{form:bracketbeta}), we get $c_{11}^{12}=0.$

Using (\ref{ec: braided equation}) and $c_{12}^{21}\neq0$, we get
$c_{11}^{22}=0.$

Using (\ref{form:antisymmetrybeta}), we land in the sixth case in
Table \ref{tab: Im(Beta)=00003D1}.

We have so proved that there exists a basis $x_{1},x_{2}$ for $V$
and $\gamma\in K$ such that $c$, $\beta$ and $f$ take one of the
canonical forms (up to an isomorphism of lifted QLie algebras) in
Table \ref{tab: Im(Beta)=00003D1}. One easily checks that the minimal
polynomials are those listed.

Let us check that the braided vector spaces in Table \ref{tab: Im(Beta)=00003D1}
are not mutual isomorphic. Clearly if the minimal polynomials are
different, the corresponding braided vector spaces can not be isomorphic.

Consider a generic basis change $\alpha:V\rightarrow V$ given by
$x_{1}=a_{1,1}x_{1}^{\prime}+a_{2,1}x_{2}^{\prime},x_{2}=a_{1,2}x_{1}^{\prime}+a_{2,2}x_{2}^{\prime}$
hence\[
\alpha=\begin{pmatrix}a_{1,1} & a_{1,2}\\
a_{2,1} & a_{2,2}\end{pmatrix}.\]
 Denote by $c^{\prime}$ and $\beta^{\prime}$ the matrices corresponding
to $c$ and $\beta$ in the new basis. Note that $c^{\prime}=\left(\alpha\otimes\alpha\right)\circ c\circ\left(\alpha^{-1}\otimes\alpha^{-1}\right)$
and $\beta^{\prime}=\alpha\circ\beta\circ\left(\alpha^{-1}\otimes\alpha^{-1}\right).$

Case A) $\beta$ is as in the fourth case in Table \ref{tab: Im(Beta)=00003D1}.

Then $\beta^{\prime}$ has zero lower row if and only if $a_{2,1}=0$.
This condition implies\[
\beta^{\prime}=\begin{pmatrix}0 & 0 & 0 & a_{1,1}/\left(a_{2,2}\right)^{2}\\
0 & 0 & 0 & 0\end{pmatrix}.\]
Hence the fourth case in Table \ref{tab: Im(Beta)=00003D1} is not
isomorphic to any other.

Case B) $\beta$ is as in the sixth case in Table \ref{tab: Im(Beta)=00003D1}.

Since $\gamma\neq1$ and $\alpha$ is invertible, then $\beta^{\prime}$
has zero lower row if and only if $a_{2,1}=0$. This condition implies\[
\beta^{\prime}=\begin{pmatrix}0 & 1/a_{2,2} & -\gamma/a_{2,2} & a_{1,2}\left(\gamma-1\right)/\left(a_{2,2}\right)^{2}\\
0 & 0 & 0 & 0\end{pmatrix}.\]
Hence, since $\gamma\neq1,$ the sixth case in Table \ref{tab: Im(Beta)=00003D1}
is not isomorphic to any other.

Case C) $\beta$ is of the form\[
\beta=\begin{pmatrix}0 & 1 & -1 & 0\\
0 & 0 & 0 & 0\end{pmatrix}.\]
Then $\beta^{\prime}$ takes the same form if and only if $a_{2,1}=0$
and $a_{2,2}=1$. Let $c$ be one of the braiding corresponding to
$\beta$ as in Table \ref{tab: Im(Beta)=00003D1}. Then $c$ and $c'$ take the form 
\[
c=\begin{pmatrix}x & y & w & z\\
0 & 0 & 1 & 0\\
0 & 1 & 0 & 0\\
0 & 0 & 0 & 1\end{pmatrix}
\quad\text{ and }\quad  
c'=
 \begin{pmatrix}
   x & y' & w' & z'\\
   0 & 0 & 1 & 0\\
   0 & 1 & 0 & 0\\
   0 & 0 & 0 & 1
 \end{pmatrix}
\]
where
\begin{equation*}
  \begin{split}
    y'&=a_{1,1}y+a_{1,2}\left(1-x\right),\\
    w'&=a_{1,1}w+a_{1,2}\left(1-x\right),\\
    z'&=\left(a_{1,1}\right)^{2}z-a_{1,1}a_{1,2}y-a_{1,1}a_{1,2}w-\left(a_{1,2}\right)^{2}\left(1-x\right).
  \end{split}
\end{equation*}

Clearly each of the first, the third, the fifth and the seventh case
in Table \ref{tab: Im(Beta)=00003D1} is only autoisomorphic.

Finally, the invariance of the minimal polynomial shows that the second
and the eighth cases in Table \ref{tab: Im(Beta)=00003D1} are only
autoisomorphic too.

Let us prove that, if $c$ and $\beta$ are as in any of the cases
of Table \ref{tab: Im(Beta)=00003D1}, then $\left(V,c,\beta\right)$
is a lifted QLie algebra. It is straightforward to check that all
of them fulfills (\ref{ec: braided equation}), (\ref{form:antisymmetrybeta})
and (\ref{form:bracketbeta}). It remains to prove that they verify
also (\ref{form:Jacobibeta}) i.e. that $\beta\left(\beta_{1}-\beta_{2}\right)\left(\overline{E_{2}\left(V,c\right)}\right)=0.$

If $c$ and $\beta$ are as in the first case of Table \ref{tab: Im(Beta)=00003D1},
then \[
E_{2}\left(V,c\right)=\mathrm{Im}\left(c-\mathrm{Id}_{V^{\otimes2}}\right)=K\left(x_{2}\otimes x_{1}-x_{1}\otimes x_{2}\right)\]
 so that $\overline{E_{2}\left(V,c\right)}=0$ whence (\ref{form:Jacobibeta})
is satisfied. We have\[
U_{Q}\left(V,c,\beta\right):=\frac{T\left(V,c\right)}{\left(\left(c-\mathrm{Id}_{V^{\otimes2}}\right)\left(z\right)-\beta\left(z\right)\mid z\in V^{\otimes2}\right)}=\frac{T\left(V,c\right)}{\left(x_{2}x_{1}-x_{1}x_{2}+x_{1}\right)}.\]

If $c$ and $\beta$ are as in the second case of Table \ref{tab: Im(Beta)=00003D1},
then \[
E_{2}\left(V,c\right)=\mathrm{Im}\left(c-\mathrm{Id}_{V^{\otimes2}}\right)=K\left(x_{1}\otimes x_{1}-x_{2}\otimes x_{1}+x_{1}\otimes x_{2}\right)\]
 so that $\overline{E_{2}\left(V,c\right)}=0$ whence (\ref{form:Jacobibeta})
is satisfied. We have\[
U_{Q}\left(V,c,\beta\right):=\frac{T\left(V,c\right)}{\left(\left(c-\mathrm{Id}_{V^{\otimes2}}\right)\left(z\right)-\beta\left(z\right)\mid z\in V^{\otimes2}\right)}=\frac{T\left(V,c\right)}{\left(\left(x_{1}\right)^{2}+x_{1}x_{2}-x_{2}x_{1}-x_{1}\right)}.\]

If $c$ and $\beta$ are as in the third case of Table \ref{tab: Im(Beta)=00003D1},
then \[
E_{2}\left(V,c\right)=\mathrm{Im}\left(c-\mathrm{Id}_{V^{\otimes2}}\right)=Kx_{1}\otimes x_{1}\oplus K\left(x_{2}\otimes x_{1}-x_{1}\otimes x_{2}\right)\]
 so that $\overline{E_{2}\left(V,c\right)}=Kx_{1}\otimes x_{1}\otimes x_{1}$.
We have that $\beta_{1}\left(\overline{E_{2}\left(V,c\right)}\right)=0=\beta_{2}\left(\overline{E_{2}\left(V,c\right)}\right),$
so that (\ref{form:Jacobibeta}) is satisfied. We have\[
U_{Q}\left(V,c,\beta\right):=\frac{T\left(V,c\right)}{\left(\left(c-\mathrm{Id}_{V^{\otimes2}}\right)\left(z\right)-\beta\left(z\right)\mid z\in V^{\otimes2}\right)}=\frac{T\left(V,c\right)}{\left(\begin{array}{c}
\left(x_{1}\right)^{2},\\
x_{2}x_{1}-x_{1}x_{2}+x_{1}\end{array}\right)}.\]

If $c$ and $\beta$ are as in the fourth case of Table \ref{tab: Im(Beta)=00003D1},
then \[
E_{2}\left(V,c\right)=\mathrm{Im}\left(c-\mathrm{Id}_{V^{\otimes2}}\right)=K\left(\gamma x_{1}\otimes x_{1}-2x_{2}\otimes x_{2}\right)\oplus K\left(x_{2}\otimes x_{1}-x_{1}\otimes x_{2}\right).\]
 We have $\beta\left(\beta_{1}-\beta_{2}\right)=0$ so that (\ref{form:Jacobibeta})
is satisfied. We have\[
U_{Q}\left(V,c,\beta\right):=\frac{T\left(V,c\right)}{\left(\left(c-\mathrm{Id}_{V^{\otimes2}}\right)\left(z\right)-\beta\left(z\right)\mid z\in V^{\otimes2}\right)}=\frac{T\left(V,c\right)}{\left(\begin{array}{c}
x_{2}x_{1}-x_{1}x_{2},\\
\gamma\left(x_{1}\right)^{2}-2\left(x_{2}\right)^{2}-x_{1}\end{array}\right)}.\]

If $c$ and $\beta$ are as in the fifth case of Table \ref{tab: Im(Beta)=00003D1},
then \[
E_{2}\left(V,c\right)=\mathrm{Im}\left(c^{2}-c\right)=K\left(x_{1}\otimes x_{1}-x_{2}\otimes x_{1}+x_{1}\otimes x_{2}\right)\]
 so that $\overline{E_{2}\left(V,c\right)}=0$ whence (\ref{form:Jacobibeta})
is satisfied. We have\[
U_{Q}\left(V,c,\beta\right):=\frac{T\left(V,c\right)}{\left(\left(c^{2}-c\right)\left(z\right)-\beta\left(z\right)\mid z\in V^{\otimes2}\right)}=\frac{T\left(V,c\right)}{\left(\left(x_{1}\right)^{2}-x_{2}x_{1}+x_{1}x_{2}+x_{1}\right)}.\]

If $c$ and $\beta$ are as in the sixth case of Table \ref{tab: Im(Beta)=00003D1},
then \[
E_{2}\left(V,c\right)=\mathrm{Im}\left(c^{2}-c\right)=K\left(\gamma x_{2}\otimes x_{1}-x_{1}\otimes x_{2}\right)\]
 so that $\overline{E_{2}\left(V,c\right)}=0$ whence (\ref{form:Jacobibeta})
is satisfied. We have\[
U_{Q}\left(V,c,\beta\right):=\frac{T\left(V,c\right)}{\left(\left(c^{2}-c\right)\left(z\right)-\beta\left(z\right)\mid z\in V^{\otimes2}\right)}=\frac{T\left(V,c\right)}{\left(-\gamma x_{2}x_{1}+x_{1}x_{2}+\gamma x_{1}\right)}.\]

If $c$ and $\beta$ are as in the seventh case of Table \ref{tab: Im(Beta)=00003D1},
then \[
E_{2}\left(V,c\right)=\mathrm{Im}\left(\left(c-\mathrm{Id}_{V^{\otimes2}}\right)\left(c-\gamma\mathrm{Id}_{V^{\otimes2}}\right)\right)=K\left(x_{2}\otimes x_{1}-x_{1}\otimes x_{2}\right)\]
 so that $\overline{E_{2}\left(V,c\right)}=0$ whence (\ref{form:Jacobibeta})
is satisfied. We have
\begin{equation*}
  \begin{split}
    U_{Q}\left(V,c,\beta\right)
     &:=\frac{T\left(V,c\right)}{\left(\left(c-\mathrm{Id}_{V^{\otimes2}}\right)\left(c-\gamma\mathrm{Id}_{V^{\otimes2}}\right)\left(z\right)-\beta\left(z\right)\mid z\in V^{\otimes2}\right)}\\
     &=\frac{T\left(V,c\right)}{\left(\left(1+\gamma\right)\left(x_{2}x_{1}-x_{1}x_{2}\right)-x_{1}\right)}.
  \end{split}
\end{equation*}

If $c$ and $\beta$ are as in the eighth case of Table \ref{tab: Im(Beta)=00003D1},
then \[
E_{2}\left(V,c\right)=\mathrm{Im}\left(\left(c-\mathrm{Id}_{V^{\otimes2}}\right)^{2}\right).\]
Now $E_{2}\left(V,c\right)=K\left(x_{2}\otimes x_{1}-x_{1}\otimes x_{2}\right)$
so that $\overline{E_{2}\left(V,c\right)}=0$ whence (\ref{form:Jacobibeta})
is satisfied. We have \begin{align*}
U_{Q}\left(V,c,\beta\right) & :=\frac{T\left(V,c\right)}{\Big(\left(c-\mathrm{Id}_{V^{\otimes2}}\right)^{2}\left(z\right)-\beta\left(z\right)\mid z\in V^{\otimes2}\Big)}\\
 & =\frac{T\left(V,c\right)}{\left(\begin{array}{c}
2\left(x_{2}x_{1}-x_{1}x_{2}\right)-x_{1}\end{array}\right)}.\end{align*}

\end{proof}
Next result is useful in conjunction with Theorem \ref{teo: magnum}.
\begin{thm}
\label{thm:tablequadr}Let $\left(V,c\right)$ be any braided vector
space as in Table \ref{tab: Im(Beta)=00003D1}. Then: \end{thm}
\begin{itemize}
\item in all cases but 7 one has $S_{Q}\left(V,c\right)=\emph{B}\left(V,c\right)$
if and only if $\mathrm{char}(K)=0$; 
\item in case 7 one has $S_{Q}\left(V,c\right)=\emph{B}\left(V,c\right)$
if and only if $\mathrm{char}(K)=0$ and $\gamma$ is not a root of
unity.\end{itemize}
\begin{proof}
Let $S_{0}:=S_{Q}(V,c_{0})$ where $c_{0}$ is the canonical flip
map and set $\Delta_{0}:=\Delta_{S_{0}}.$ Set $S:=S_{Q}\left(V,c\right)$.
First note that $S=\emph{B}\left(V,c\right)$ if and only if $P\left(S\right)=Kx_{1}+Kx_{2}$.

Suppose $\mathrm{char}(K)=0.$ Then, by \cite[Theorem 2.17]{AMS-Hecke-Type},
for $\left(V,c\right)$ as in the cases 1,2,3 and 4 of Table \ref{tab: Im(Beta)=00003D1},
we get $S=\emph{B}\left(V,c\right)$. 

Let us concern the cases 5, 6 and 7 of Table \ref{tab: Im(Beta)=00003D1}.
Then $c$ has matrix of the form \[
\begin{pmatrix}q_{1} & q & 0 & 0\\
0 & 0 & q_{12} & 0\\
0 & q_{21} & 0 & 0\\
0 & 0 & 0 & q_{2}\end{pmatrix}\]
where $q_{21}=q_{12}^{-1}$ and $q_{2}=1$. Then \[
\Delta\left(x_{i}^{n}\right)=\sum_{0\leq t\leq n}\binom{n}{t}_{q_{i}}\left(x_{i}^{t}\otimes x_{i}^{n-t}\right)\]
so that \begin{eqnarray*}
 &  & \Delta\left(x_{1}^{n_{1}}x_{2}^{n_{2}}\right)\\
 & = & \Delta\left(x_{1}^{n_{1}}\right)\Delta\left(x_{2}^{n_{2}}\right)\\
 & = & \sum_{\begin{array}{c}
0\leq t_{1}\leq n_{1}\\
0\leq t_{2}\leq n_{2}\end{array}}\binom{n_{1}}{t_{1}}_{q_{1}}\binom{n_{2}}{t_{2}}_{q_{2}}\left(x_{1}^{t_{1}}\otimes x_{1}^{n_{1}-t_{1}}\right)\left(x_{2}^{t_{2}}\otimes x_{2}^{n_{2}-t_{2}}\right)\\
 & = & \sum_{\begin{array}{c}
0\leq t_{1}\leq n_{1}\\
0\leq t_{2}\leq n_{2}\end{array}}\binom{n_{1}}{t_{1}}_{q_{1}}\binom{n_{2}}{t_{2}}_{q_{2}}x_{1}^{t_{1}}c\left(x_{1}^{n_{1}-t_{1}}\otimes x_{2}^{t_{2}}\right)x_{2}^{n_{2}-t_{2}}\\
 & = & \sum_{\begin{array}{c}
0\leq t_{1}\leq n_{1}\\
0\leq t_{2}\leq n_{2}\end{array}}\binom{n_{1}}{t_{1}}_{q_{1}}\binom{n_{2}}{t_{2}}_{q_{2}}q_{12}^{\left(n_{1}-t_{1}\right)t_{2}}x_{1}^{t_{1}}\left(x_{2}^{t_{2}}\otimes x_{1}^{n_{1}-t_{1}}\right)x_{2}^{n_{2}-t_{2}}\\
 & = & \sum_{\begin{array}{c}
0\leq t_{1}\leq n_{1}\\
0\leq t_{2}\leq n_{2}\end{array}}\binom{n_{1}}{t_{1}}_{q_{1}}\binom{n_{2}}{t_{2}}_{q_{2}}q_{12}^{\left(n_{1}-t_{1}\right)t_{2}}x_{1}^{t_{1}}x_{2}^{t_{2}}\otimes x_{1}^{n_{1}-t_{1}}x_{2}^{n_{2}-t_{2}}\end{eqnarray*}

so that \[
\Delta\left(x_{1}^{n_{1}}x_{2}^{n_{2}}\right)=\sum_{\begin{array}{c}
0\leq t_{1}\leq n_{1}\\
0\leq t_{2}\leq n_{2}\end{array}}\binom{n_{1}}{t_{1}}_{q_{1}}\binom{n_{2}}{t_{2}}_{q_{2}}q_{12}^{\left(n_{1}-t_{1}\right)t_{2}}x_{1}^{t_{1}}x_{2}^{t_{2}}\otimes x_{1}^{n_{1}-t_{1}}x_{2}^{n_{2}-t_{2}}.\]
Let $n>1$ and $z\in P\left(S\right)\cap S_{n}$. Then \[
z=\sum_{\begin{array}{c}
0\leq n_{1},n_{2}\\
n_{1}+n_{2}=n\end{array}}\lambda_{n_{1},n_{2}}x_{1}^{n_{1}}x_{2}^{n_{2}}.\]
Then $z\otimes1+1\otimes z=\Delta\left(z\right)$ rewrites as
\begin{eqnarray*}
 &  & \sum_{\begin{array}{c}
            0\leq n_{1},n_{2}\\
            n_{1}+n_{2}=n
            \end{array}}
          \lambda_{n_{1},n_{2}}x_{1}^{n_{1}}x_{2}^{n_{2}}\otimes1
      +1\otimes\sum_{\begin{array}{c}
                       0\leq n_{1},n_{2}\\ 
                       n_{1}+n_{2}=n
                     \end{array}}
                \lambda_{n_{1},n_{2}}x_{1}^{n_{1}}x_{2}^{n_{2}}\\
 & = & \sum_{\begin{array}{c}
                0\leq n_{1},n_{2}\\ 
                n_{1}+n_{2}=n\\
                0\leq t_{1}\leq n_{1}\\
                0\leq t_{2}\leq n_{2}
              \end{array}}
            \lambda_{n_{1},n_{2}}\binom{n_{1}}{t_{1}}_{q_{1}}\binom{n_{2}}{t_{2}}_{q_{2}}q_{12}^{\left(n_{1}-t_{1}\right)t_{2}}x_{1}^{t_{1}}x_{2}^{t_{2}}\otimes x_{1}^{n_{1}-t_{1}}x_{2}^{n_{2}-t_{2}}.
\end{eqnarray*}
Since $n>1$ we can always find a couple $\left(t_{1},t_{2}\right)\neq\left(0,0\right),\left(n_{1},n_{2}\right)$
so that we get \[
\lambda_{n_{1},n_{2}}\binom{n_{1}}{t_{1}}_{q_{1}}\binom{n_{2}}{t_{2}}_{q_{2}}q_{12}^{\left(n_{1}-t_{1}\right)t_{2}}=0.\]
Since $\mathrm{char}(K)=0$ and in case 7 one has that $\gamma$ is
not root of unity, we get $\lambda_{n_{1},n_{2}}=0$. Hence $z=0$.
We have so proved that $P\left(S\right)\cap S_{n}=\left\{ 0\right\} $
for all $n>1.$ Hence $P\left(S\right)=Kx_{1}+Kx_{2}$ so that $S=B\left(V,c\right)$.

Let us consider the eighth case of Table \ref{tab: Im(Beta)=00003D1}.
Let us prove inductively that \begin{equation}
c\left(x_{2}\otimes x_{2}^{n}\right)=x_{2}^{n}\otimes x_{2}+n\gamma x_{1}x_{2}^{n-1}\otimes x_{1}.\label{eq:cx2}\end{equation}
For $n=0$ there is nothing to prove. For $n>0$ we have\begin{eqnarray*}
c\left(x_{2}\otimes x_{2}^{n+1}\right) & = & c\left(x_{2}\otimes x_{2}x_{2}^{n}\right)\\
 & = & \left(m\otimes S\right)\left(S\otimes c\right)\left(c\otimes S\right)\left(x_{2}\otimes x_{2}\otimes x_{2}^{n}\right)\\
 & = & x_{2}c\left(x_{2}\otimes x_{2}^{n}\right)+\gamma x_{1}c\left(x_{1}\otimes x_{2}^{n}\right)\\
 & = & x_{2}^{n+1}\otimes x_{2}+n\gamma x_{2}x_{1}x_{2}^{n-1}\otimes x_{1}+\gamma x_{1}x_{2}^{n}\otimes x_{1}\\
 & = & x_{2}^{n+1}\otimes x_{2}+\left(n+1\right)\gamma x_{1}x_{2}^{n}\otimes x_{1}.\end{eqnarray*}

For $n\geq0$, let us prove there exist $\alpha_{t}(n)\in K$ such
that\[
\Delta\left(x_{2}^{n}\right)=\sum_{0\leq t}\alpha_{t}(n)\gamma^{t}\left(x_{1}^{t}\otimes x_{1}^{t}\right)\Delta_{0}\left(x_{2}^{n-2t}\right)\]
where we assume $x_{2}^{t}:=0$ for $t<0$. 

For $n=0$ we get $\alpha_{t}(0)=\delta_{t,0}$. For $n=1$, one easily
get $\alpha_{t}(1)=\delta_{t,0}$. Let $n>1$ and assume the statement is
true for $n$. First we compute
\begin{eqnarray*}
  \left(1\otimes x_{2}\right)\Delta_{0}\left(x_{2}^{m}\right) & = & \sum_{0\leq i\leq m}\binom{m}{i}\left(1\otimes x_{2}\right)\left(x_{2}^{i}\otimes x_{2}^{m-i}\right)\\
  & = & \sum_{0\leq i\leq m}\binom{m}{i}c\left(x_{2}\otimes x_{2}^{i}\right)x_{2}^{m-i}\\
  & \overset{\eqref{eq:cx2}}{=} & \sum_{0\leq i\leq m}\binom{m}{i}\left(x_{2}^{i}\otimes x_{2}^{m+1-i}+i\gamma x_{1}x_{2}^{i-1}\otimes x_{1}x_{2}^{m-i}\right)\\
  & = & \sum_{0\leq i\leq m}\binom{m}{i}x_{2}^{i}\otimes x_{2}^{m+1-i}+\\
  & & +\gamma\sum_{0\leq i\leq m}\binom{m}{i}ix_{1}x_{2}^{i-1}\otimes x_{1}x_{2}^{m-i}\\
  & = & \left(1\otimes x_{2}\right)\cdot_{S_{0}\otimes S_{0}}\Delta_{0}\left(x_{2}^{m}\right)+\\
  & & +\gamma\left(x_{1}\otimes x_{1}\right)\sum_{1\leq i\leq m}\binom{m}{i}ix_{2}^{i-1}\otimes x_{2}^{m-i}\\
  & = & \left(1\otimes x_{2}\right)\cdot_{S_{0}\otimes S_{0}}\Delta_{0}\left(x_{2}^{m}\right)+m\gamma\left(x_{1}\otimes x_{1}\right)\Delta_{0}\left(x_{2}^{m-1}\right)
\end{eqnarray*}
so that \begin{equation}
\left(1\otimes x_{2}\right)\Delta_{0}\left(x_{2}^{m}\right)=\left(1\otimes x_{2}\right)\cdot_{S_{0}\otimes S_{0}}\Delta_{0}\left(x_{2}^{m}\right)+m\gamma\left(x_{1}\otimes x_{1}\right)\Delta_{0}\left(x_{2}^{m-1}\right).\label{eq:deltx2}\end{equation}
We have\begin{eqnarray*}
 &  & \Delta x_{2}^{n+1}\\
 & = & \left(\Delta x_{2}\right)\left(\Delta x_{2}^{n}\right)\\
 & = & \left(x_{2}\otimes1+1\otimes x_{2}\right)\left(\sum_{0\leq t}\alpha_{t}(n)\gamma^{t}\left(x_{1}^{t}\otimes x_{1}^{t}\right)\Delta_{0}\left(x_{2}^{n-2t}\right)\right)\\
 & = & \left(\begin{array}{c}
\sum_{0\leq t}\alpha_{t}(n)\gamma^{t}\left(x_{2}\otimes1\right)\left(x_{1}^{t}\otimes x_{1}^{t}\right)\Delta_{0}\left(x_{2}^{n-2t}\right)\\
+\sum_{0\leq t}\alpha_{t}(n)\gamma^{t}\left(1\otimes x_{2}\right)\left(x_{1}^{t}\otimes x_{1}^{t}\right)\Delta_{0}\left(x_{2}^{n-2t}\right)\end{array}\right)\\
 & = & \left(\begin{array}{c}
\sum_{0\leq t}\alpha_{t}(n)\gamma^{t}\left(x_{1}^{t}\otimes x_{1}^{t}\right)\left[\left(x_{2}\otimes1\right)\Delta_{0}\left(x_{2}^{n-2t}\right)\right]\\
+\sum_{0\leq t}\alpha_{t}(n)\gamma^{t}\left(x_{1}^{t}\otimes x_{1}^{t}\right)\left[\left(1\otimes x_{2}\right)\Delta_{0}\left(x_{2}^{n-2t}\right)\right]\end{array}\right)\\
 & \overset{\eqref{eq:deltx2}}{=} & \left(\begin{array}{c}
\sum_{0\leq t}\alpha_{t}(n)\gamma^{t}\left(x_{1}^{t}\otimes x_{1}^{t}\right)\left[\left(x_{2}\otimes1\right)\cdot_{S_{0}\otimes S_{0}}\Delta_{0}\left(x_{2}^{n-2t}\right)\right]\\
+\sum_{0\leq t}\alpha_{t}(n)\gamma^{t}\left(x_{1}^{t}\otimes x_{1}^{t}\right)\left[\left(1\otimes x_{2}\right)\cdot_{S_{0}\otimes S_{0}}\Delta_{0}\left(x_{2}^{n-2t}\right)\right]\\
+\sum_{0\leq t}\alpha_{t}(n)\gamma^{t}\left(x_{1}^{t}\otimes x_{1}^{t}\right)\left[\left(n-2t\right)\gamma\left(x_{1}\otimes x_{1}\right)\Delta_{0}\left(x_{2}^{n-2t-1}\right)\right]\end{array}\right)\\
 & = & \left(\begin{array}{c}
\sum_{0\leq t}\alpha_{t}(n)\gamma^{t}\left(x_{1}^{t}\otimes x_{1}^{t}\right)\left[\left(x_{2}\otimes1+1\otimes x_{2}\right)\cdot_{S_{0}\otimes S_{0}}\Delta_{0}\left(x_{2}^{n-2t}\right)\right]\\
+\sum_{0\leq t}\alpha_{t}(n)\left(n-2t\right)\gamma^{t+1}\left(x_{1}^{t+1}\otimes x_{1}^{t+1}\right)\Delta_{0}\left(x_{2}^{n-2t-1}\right)\end{array}\right)\\
 & = & \left(\begin{array}{c}
\sum_{0\leq t}\alpha_{t}(n)\gamma^{t}\left(x_{1}^{t}\otimes x_{1}^{t}\right)\Delta_{0}\left(x_{2}^{n+1-2t}\right)\\
+\sum_{0\leq t}\alpha_{t}(n)\left(n-2t\right)\gamma^{t+1}\left(x_{1}^{t+1}\otimes x_{1}^{t+1}\right)\Delta_{0}\left(x_{2}^{n-2t-1}\right)\end{array}\right)\\
 & = & \left(\begin{array}{c}
\sum_{0\leq t}\alpha_{t}(n)\gamma^{t}\left(x_{1}^{t}\otimes x_{1}^{t}\right)\Delta_{0}\left(x_{2}^{n+1-2t}\right)\\
+\sum_{1\leq s}\alpha_{s-1}(n)\left(n+2-2s\right)\gamma^{s}\left(x_{1}^{s}\otimes x_{1}^{s}\right)\Delta_{0}\left(x_{2}^{n+1-2s}\right)\end{array}\right)\\
 & = & \begin{array}{c}
\sum_{0\leq t}\alpha_{t}(n+1)\gamma^{t}\left(x_{1}^{t}\otimes x_{1}^{t}\right)\Delta_{0}\left(x_{2}^{n+1-2t}\right)\end{array}\end{eqnarray*}
where \[
\alpha_{t}(n+1)=\begin{cases}
\alpha_{t}(n) & \qquad\text{if }t=0,\\
\alpha_{t}(n)+\alpha_{t-1}(n)\left(n+2-2t\right) & \qquad\text{if }t>0.\end{cases}\]
Since $\alpha_{0}(0)=1$ we get $\alpha_{0}(n)=1$ for all $n\geq0$
and, since for all $t>0$, $\alpha_{t}(1)=0$, one easily get $\alpha_{t}(n)=0$
whenever $n<2t.$ 

Now we have\begin{eqnarray*}
\Delta\left(x_{1}^{n_{1}}x_{2}^{n_{2}}\right) & = & \Delta\left(x_{1}^{n_{1}}\right)\Delta\left(x_{2}^{n_{2}}\right)\\
 & = & \Delta\left(x_{1}^{n_{1}}\right)\left[\sum_{0\leq t}\alpha_{t}(n_{2})\gamma^{t}\left(x_{1}^{t}\otimes x_{1}^{t}\right)\Delta_{0}\left(x_{2}^{n_{2}-2t}\right)\right]\\
 & = & \sum_{0\leq t}\alpha_{t}(n_{2})\gamma^{t}\left(x_{1}^{t}\otimes x_{1}^{t}\right)\Delta\left(x_{1}^{n_{1}}\right)\Delta_{0}\left(x_{2}^{n_{2}-2t}\right)\\
 & = & \sum_{0\leq t}\alpha_{t}(n_{2})\gamma^{t}\left(x_{1}^{t}\otimes x_{1}^{t}\right)\left[\Delta_{0}\left(x_{1}^{n_{1}}\right)\cdot_{S_{0}\otimes S_{0}}\Delta_{0}\left(x_{2}^{n_{2}-2t}\right)\right]\\
 & = & \sum_{0\leq t}\alpha_{t}(n_{2})\gamma^{t}\left(x_{1}^{t}\otimes x_{1}^{t}\right)\Delta_{0}\left(x_{1}^{n_{1}}x_{2}^{n_{2}-2t}\right).\end{eqnarray*}
Let $n>1$ and $z\in P\left(S\right)\cap S_{n}$. Then \[
z=\sum_{\begin{array}{c}
0\leq n_{1},n_{2}\\
n_{1}+n_{2}=n\end{array}}\lambda_{n_{1},n_{2}}x_{1}^{n_{1}}x_{2}^{n_{2}}.\]
Then $z\otimes1+1\otimes z=\Delta\left(z\right)$ rewrites as\begin{eqnarray*}
 &  & \sum_{\begin{array}{c}
0\leq n_{1},n_{2}\\
n_{1}+n_{2}=n\end{array}}\lambda_{n_{1},n_{2}}x_{1}^{n_{1}}x_{2}^{n_{2}}\otimes1+1\otimes\sum_{\begin{array}{c}
0\leq n_{1},n_{2}\\
n_{1}+n_{2}=n\end{array}}\lambda_{n_{1},n_{2}}x_{1}^{n_{1}}x_{2}^{n_{2}}\\
 & = & \sum_{\begin{array}{c}
0\leq n_{1},n_{2}\\
n_{1}+n_{2}=n\\
0\leq t\end{array}}\lambda_{n_{1},n_{2}}\alpha_{t}(n_{2})\gamma^{t}\left(x_{1}^{t}\otimes x_{1}^{t}\right)\Delta_{0}\left(x_{1}^{n_{1}}x_{2}^{n_{2}-2t}\right).\end{eqnarray*}
Note that 
  \begin{multline}
    \left(x_{1}^{t}\otimes x_{1}^{t}\right)\Delta_{0}\left(x_{1}^{n_{1}}x_{2}^{n_{2}-2t}\right)=\\
    \sum_{\begin{array}{c}
      0\leq i_{1}\leq n_{1}\\
      0\leq i_{2}\leq n_{2}
    \end{array}}
    \binom{n_{1}}{i_{1}}\binom{n_{2}}{i_{2}}\left(x_{1}^{t+i_{1}}x_{2}^{i_{2}}\otimes x_{1}^{t+n_{1}-i_{1}}x_{2}^{n_{2}-i_{2}}\right).
  \end{multline}
Since $n>1$ we can always find a couple $\left(i_{1},i_{2}\right)\neq\left(0,0\right),\left(n_{1},n_{2}\right)$
so that we get \[
\lambda_{n_{1},n_{2}}\alpha_{t}(n_{2})\gamma^{t}\binom{n_{1}}{i_{1}}\binom{n_{2}}{i_{2}}=0.\]
For $t=0$, since $\mathrm{char}(K)=0$, we get $\lambda_{n_{1},n_{2}}=0$.
Hence $z=0$. We have so proved that $P\left(S\right)\cap S_{n}=\left\{ 0\right\} $
for all $n>1.$ Hence \[
P\left(S\right)=Kx_{1}+Kx_{2}\]
so that $S=B\left(V,c\right)$.

Conversely, suppose $S=\emph{B}\left(V,c\right)$. Then $P\left(S\right)=Kx_{1}+Kx_{2}$.
By contradiction, suppose $\mathrm{char}(K)=p$ for some prime $p>2$.
In the cases 1, 2, 4 and 8 of Table \ref{tab: Im(Beta)=00003D1},
we get that $x_{1}^{p}\in P\left(S\right)$ , a contradiction. In
the remaining cases we get that $x_{2}^{p}\in P\left(S\right)$, a
contradiction. Note that in case 3 the proof uses that \[
\Delta\left(x_{2}^{n}\right)=\Delta_{0}\left(x_{2}^{n}\right)+\binom{n}{2}\gamma\left(x_{1}\otimes x_{1}\right)\Delta_{0}\left(x_{2}^{n-2}\right)\]
that can be proved in a similar way as in the other implication. 

In case 7, if $\mathrm{char}(K)=0$ and $\gamma$ is a primitive $t$-root
of unity for some $t>0$, then $t>2$ (as $\gamma\neq\pm1$) and we
get $x_{1}^{t}\in P\left(S\right)$, a contradiction .\end{proof}
\begin{cor}
Let $\left(V,c\right)$ be any braided vector space as in Table \ref{tab: Im(Beta)=00003D1}.
Then the equivalent conditions of Theorem \ref{teo: magnum} are fulfilled whenever $S_{Q}\left(V,c\right)=\emph{B}\left(V,c\right)$.\end{cor}
\begin{proof}
It is enough to observe that condition (ii) of Theorem \ref{teo: magnum}
holds for all $\left(V,c\right)$ as in Table \ref{tab: Im(Beta)=00003D1}.\end{proof}
\begin{thm}
\label{thm:sunto}Let $A$ be a primitively generated connected braided
bialgebra and let $(P,c)$ be the braided vector space of primitive
elements of $A$. Assume the Nichols algebra $\emph{B}\left(P,c\right)$
is a quadratic algebra and $c$ has minimal polynomial $f\in K\left[X\right]$ of the form $f=\left(X+1\right)h$ for some $h\in K\left[X\right]$ such that
$h\left(-1\right)\neq0.$ Set $\beta:=m_Ah(c)$. Then $\left(P,c,\beta\right)$ is a lifted QLie algebra  and $A$ is isomorphic to the universal enveloping algebra $U_{Q}\left(V,c,\beta\right)$.
Moreover assume 
\begin{enumerate}
\item [1)] $\mathrm{char}\left(K\right)\neq2$;
\item [2)] $P$ is two dimensional;
\item [3)] $\dim\left(\mathrm{Im}\beta\right)=1.$ 
\end{enumerate}
Then, there exists a basis $x_{1},x_{2}$ for $P$ and $\gamma\in K$
such that $c$, $\beta,$ $U_{Q}\left(P,c,\beta\right)$ and $f$
take only one of the canonical forms (up to an isomorphism of lifted
QLie algebras) in Table \ref{tab: Im(Beta)=00003D1}.\end{thm}
\begin{proof}
The first part follows by Theorem \ref{thm:quadratic} and the second
part by Theorem \ref{thm:teoClassif}. 
\end{proof}
\appendix

\section{Further results\label{sec:appendix}}

Although we don't have a complete classification in case $\dim\left(\mathrm{Im}\beta\right)=2$
yet, in this section we include a partial result that can help in
this direction.
\begin{thm}
If $\dim\left(\mathrm{Im}\beta\right)=2,$ then $\dim\mathrm{Im}\left(c+\mathrm{Id}_{V^{\otimes2}}\right)=2.$ \end{thm}
\begin{proof}
Assume $\dim\left(\mathrm{Im}\beta\right)=2.$ Then $\beta$ is surjective
so that $\dim\mathrm{Im}\left(h\left(c\right)\right)\geq2.$

Since $h(-1)\neq0$, by Lemma \ref{lem: utile}, $\mathrm{Im}\left(c+\mathrm{Id}_{V^{\otimes2}}\right)\oplus\mathrm{Im}\left(h\left(c\right)\right)=V^{\otimes2}$
in view of the standing hypotheses in Notation \ref{not: dim 2}.
Since $V$ is two-dimensional, we have three cases namely $\dim\mathrm{Im}\left(c+\mathrm{Id}_{V^{\otimes2}}\right)=0,1,2.$
The first case does not occur. In fact, $\dim\mathrm{Im}\left(c+\mathrm{Id}_{V^{\otimes2}}\right)=0$
implies $c=-\mathrm{Id}_{V^{\otimes2}}.$ By (\ref{form:bracketbeta}),
we get $\beta_{1}=-\beta_{2}\ $which entails $\beta=0,$ a contradiction.

Let $\delta:\mathrm{Im}\left(c+\mathrm{Id}_{V^{\otimes2}}\right)\rightarrow\mathrm{Im}\left(c+\mathrm{Id}_{V^{\otimes2}}\right)$
be the restriction of $c+\mathrm{Id}_{V^{\otimes2}}.$ Let $p_{\delta}\in K\left[X\right]$
be the minimal polynomial of $\delta.$ Then $f\mid p_{\delta}\cdot\left(X+1\right).$

$\dim\mathrm{Im}\left(c+\mathrm{Id}_{V^{\otimes2}}\right)=1)$ In
this case $\deg\left(p_{\delta}\right)=1$ so that $p_{\delta}=X-q$
for some $q\in K$ and $f=\left(X-q\right)\left(X+1\right).$ 

We need the following lemma.

\begin{lem}
\label{lem: ker(beta)}$\ker\left(\beta\right)=\mathrm{Im}\left(c+\mathrm{Id}_{V^{\otimes2}}\right)\oplus\left[\ker\left(\beta\right)\cap\mathrm{Im}\left(h\left(c\right)\right)\right].$ \end{lem}
\begin{proof}
Since $\dim\mathrm{Im}\left(c+\mathrm{Id}_{V^{\otimes2}}\right)=1$
we get $\mathrm{Im}\left(c+\mathrm{Id}_{V^{\otimes2}}\right)=Kv_{1}$
for some $v_{1}\in V^{\otimes2}\backslash\left\{ 0\right\} .$ By
(\ref{form:antisymmetrybeta}), $v_{1}\in\ker\left(\beta\right).$
Since $\beta$ is surjective, one has \[
V=\mathrm{Im}\left(\beta\right)\simeq V^{\otimes2}/\ker\left(\beta\right)\]
 so that $\dim\ker\left(\beta\right)=2.$ Hence we can complete $v_{1}$
to a basis $v_{1},v_{2}^{\prime}$ of $\ker\left(\beta\right)$. Now
$v_{2}^{\prime}\in V^{\otimes2}=\mathrm{Im}\left(c+\mathrm{Id}_{V^{\otimes2}}\right)\oplus\mathrm{Im}\left(h\left(c\right)\right)$
so that $v_{2}^{\prime}=av_{1}+v_{2}$ for some $v_{2}\in\mathrm{Im}\left(h\left(c\right)\right)$.
Moreover $v_{2}=-av_{1}+v_{2}^{\prime}\in\ker\left(\beta\right)$
whence $v_{2}\in\ker\left(\beta\right)\cap\mathrm{Im}\left(h\left(c\right)\right)$.
Clearly $v_{1},v_{2}$ is a basis for $\ker\left(\beta\right)$. 
\end{proof}
In view of Lemma \ref{lem: ker(beta)}, we have that $\ker\left(\beta\right)$
has a basis $v_{1},v_{2}$ where $v_{1}\in\mathrm{Im}\left(c+\mathrm{Id}_{V^{\otimes2}}\right)=\mathrm{\ker}\left(h\left(c\right)\right)=\mathrm{\ker}\left(c-q\mathrm{Id}_{V^{\otimes2}}\right)$
and $v_{2}\in\left[\ker\left(\beta\right)\cap\mathrm{Im}\left(h\left(c\right)\right)\right]$.
Complete $v_{2}$ to a basis $v_{2},v_{3},v_{4}$ of $\mathrm{Im}\left(h\left(c\right)\right)=\ker\left(c+\mathrm{Id}_{V^{\otimes2}}\right)$.
Clearly $v_{1},v_{2},v_{3},v_{4}$ is a basis of $V^{\otimes2}$.
We set $x_{1}:=\beta\left(v_{3}\right)$ and $x_{2}:=\beta\left(v_{4}\right)$.
Obviously $x_{1},x_{2}$ is a basis of $V$.

Suppose $x_{1}\otimes x_{1}\notin\mathrm{Im}\left(h\left(c\right)\right).$
Then $x_{1}\otimes x_{1},v_{2},v_{3},v_{4}$ is a basis of $V^{\otimes2}.$
Thus for each $i,j\in\left\{ 1,2\right\} $ there are $\xi_{t}^{i,j}\in K$
such that\[
x_{i}\otimes x_{j}=\xi_{1}^{i,j}x_{1}\otimes x_{1}+\xi_{2}^{i,j}v_{2}+\xi_{3}^{i,j}v_{3}+\xi_{4}^{i,j}v_{4}.\]
Clearly $\xi_{t}^{1,1}=\delta_{1,t}.$ Now\begin{gather*}
\left(c+\mathrm{Id}_{V^{\otimes2}}\right)\left(x_{i}\otimes x_{j}\right)=\xi_{1}^{i,j}\left(c+\mathrm{Id}_{V^{\otimes2}}\right)\left(x_{1}\otimes x_{1}\right),\\
\beta\left(x_{i}\otimes x_{j}\right)=\xi_{1}^{i,j}\beta\left(x_{1}\otimes x_{1}\right)+\xi_{3}^{i,j}x_{1}+\xi_{4}^{i,j}x_{2}.\end{gather*}
By setting $\left(c+\mathrm{Id}_{V^{\otimes2}}\right)\left(x_{1}\otimes x_{1}\right)=\sum_{m,n}v^{m,n}x_{m}\otimes x_{n}$
and $\beta\left(x_{1}\otimes x_{1}\right)=\sum_{m}s^{m}x_{m}$ we
get\begin{eqnarray*}
c\left(x_{i}\otimes x_{j}\right) & = & -x_{i}\otimes x_{j}+\sum_{m,n}\xi_{1}^{i,j}v^{m,n}x_{m}\otimes x_{n},\\
\beta\left(x_{i}\otimes x_{j}\right) & = & \sum_{m}\xi_{1}^{i,j}s^{m}x_{m}+\xi_{3}^{i,j}x_{1}+\xi_{4}^{i,j}x_{2}.\end{eqnarray*}
In matrix form we have
\begin{multline}
  c=
  \begin{pmatrix}
    -1 & 0 & 0 & 0\\
    0 & -1 & 0 & 0\\
    0 & 0 & -1 & 0\\
    0 & 0 & 0 & -1
  \end{pmatrix}
  +\\
  \begin{pmatrix}
    v^{1,1} & 0 & 0 & 0\\
    0 & v^{2,1} & 0 & 0\\
    0 & 0 & v^{1,2} & 0\\
    0 & 0 & 0 & v^{2,2}
  \end{pmatrix}
  \begin{pmatrix}
    1 & 1 & 1 & 1\\
    1 & 1 & 1 & 1\\
    1 & 1 & 1 & 1\\
    1 & 1 & 1 & 1
  \end{pmatrix}
  \begin{pmatrix}
    1 & 0 & 0 & 0\\
    0 & \xi_{1}^{2,1} & 0 & 0\\
    0 & 0 & \xi_{1}^{1,2} & 0\\
    0 & 0 & 0 & \xi_{1}^{2,2}
  \end{pmatrix},
\end{multline}
and
\begin{equation*}
  \beta=
   \begin{pmatrix}
     0 & \xi_{3}^{2,1} & \xi_{3}^{1,2} & \xi_{3}^{2,2}\\
     0 & \xi_{4}^{2,1} & \xi_{4}^{1,2} & \xi_{4}^{2,2}
   \end{pmatrix}
   +
   \begin{pmatrix}
     s^{1} & 0\\
     0 & s^{2}
   \end{pmatrix}
   \begin{pmatrix}
     1 & 1 & 1 & 1\\
     1 & 1 & 1 & 1
   \end{pmatrix}
   \begin{pmatrix}
     1 & 0 & 0 & 0\\
     0 & \xi_{1}^{2,1} & 0 & 0\\
     0 & 0 & \xi_{1}^{1,2} & 0\\
     0 & 0 & 0 & \xi_{1}^{2,2}
   \end{pmatrix}.
\end{equation*}

We set \begin{align*}
\mathbf{1} & :=\text{identity $4\times4$ matrix},\\
\mathcal{V} & :=\text{diag}\{v^{1,1},v^{2,1},v^{1,2},v^{2,2}\},\\
\Xi & :=\text{diag}\{1,\xi_{1}^{2,1},\xi_{1}^{1,2},\xi_{1}^{2,2}\}\text{ and}\\
\mathcal{U} & :=4\times4\text{ matrix with 1 in each entry}\end{align*}
 moreover, we set $\mathcal{M}:=\mathcal{VU}\Xi$.

Now 
\begin{equation*}
  \begin{split}
     f\left(c\right)=0  &\Leftrightarrow   \left(c-q\mathrm{Id}_{V^{\otimes2}}\right)\left(c+\mathrm{Id}_{V^{\otimes2}}\right)=0\\ 
    &\Leftrightarrow   \mathcal{M}\left[\mathcal{M}-\left(q+1\right)\mathbf{1}\right]=0 \\
    &\Leftrightarrow   \mathcal{VU}\Xi\left[\mathcal{VU}\Xi-\left(q+1\right)\mathbf{1}\right]=0 \\
    &\Leftrightarrow   \mathcal{V}\big[\mathcal{U}\Xi\mathcal{VU}-\left(q+1\right)\mathcal{U}\big]\Xi=0.
  \end{split}
\end{equation*}

It is straightforward to check that \begin{equation}
\mathcal{U}D\mathcal{U}=Tr\left(D\right)\cdot\mathcal{U},\text{ for each diagonal matrix }D.\text{ }\label{form: diag Stumbo}\end{equation}

By (\ref{form: diag Stumbo}), we have $\mathcal{U}\Xi\mathcal{VU}-\left(q+1\right)\mathcal{U}=x\mathcal{U}$
where $x:=Tr\left(\Xi\mathcal{V}\right)-\left(q+1\right).$ Hence
\[
\mathcal{V}\big[\mathcal{U}\Xi\mathcal{VU}-\left(q+1\right)\mathcal{U}\big]\Xi=x\mathcal{VU}\Xi=x\mathcal{M}.\]
 Since $c+\mathrm{Id}_{V^{\otimes2}}\neq0,$ we have that $\mathcal{M}$
is not zero so that $\mathcal{V}\big[\mathcal{U}\Xi\mathcal{VU}-\left(q+1\right)\mathcal{U}\big]\Xi=0$
is equivalent to $x=0,$ i.e.\ to\[
v^{11}=q+1-\left(v^{2,1}\xi_{1}^{2,1}+v^{1,2}\xi_{1}^{1,2}+v^{2,2}\xi_{1}^{2,2}\right).\]
Hence we can apply this substitution. Then one easily gets that (\ref{form:antisymmetrybeta})
is equivalent to the following extra conditions\begin{eqnarray*}
s^{1} & = & -\frac{1}{q+1}\left(v^{2,1}\xi_{3}^{2,1}+v^{1,2}\xi_{3}^{1,2}+v^{2,2}\xi_{3}^{2,2}\right),\\
s^{2} & = & -\frac{1}{q+1}\left(v^{2,1}\xi_{4}^{2,1}+v^{1,2}\xi_{4}^{1,2}+v^{2,2}\xi_{4}^{2,2}\right).\end{eqnarray*}
These conditions force $\beta=0,$ which is a contradiction (we performed
the computation by means of \cite{Axiom}). 

Clearly also assuming $x_{2}\otimes x_{2}\notin\mathrm{Im}\left(h\left(c\right)\right)$
one gets a contradiction. In conclusion we can assume $x_{i}\otimes x_{i}\in\mathrm{Im}\left(h\left(c\right)\right)$
for each $i\in\left\{ 1,2\right\} .$ Then $c$ and $\beta$ take
the form \[
c=\begin{pmatrix}-1 & c_{11}^{21} & c_{11}^{12} & 0\\
0 & c_{21}^{21} & c_{21}^{12} & 0\\
0 & c_{12}^{21} & c_{12}^{12} & 0\\
0 & c_{22}^{21} & c_{22}^{12} & -1\end{pmatrix}\qquad\text{and}\qquad\beta=\begin{pmatrix}\beta_{1}^{11} & \beta_{1}^{21} & \beta_{1}^{12} & \beta_{1}^{22}\\
\beta_{2}^{11} & \beta_{2}^{21} & \beta_{2}^{12} & \beta_{2}^{22}\end{pmatrix}.\]
We distinguish between different cases.

CASE 1) $\beta_{2}^{11}\neq0.$ In this case, by Proposition \ref{pro: normalization of beta},
we can assume $\beta_{2}^{11}=1.$

Using the left-hand side of (\ref{form:bracketbeta}), we get $c_{21}^{21}=0,c_{12}^{21}=1,c_{22}^{21}=0.$

Using the right-hand side of (\ref{form:bracketbeta}), we get a contradiction.

CASE 2) $\beta_{2}^{11}=0.$

Using the left-hand side of (\ref{form:bracketbeta}), we get $\beta_{1}^{11}=0.$

CASE 2.1) $\beta_{1}^{22}\neq0.$ In this case, by Proposition \ref{pro: normalization of beta},
we can assume $\beta_{1}^{22}=1.$

Using the left-hand side of (\ref{form:bracketbeta}), we get $c_{12}^{12}=0,$
$c_{11}^{12}=0,c_{21}^{12}=0.$

Using the right-hand side of (\ref{form:bracketbeta}), we get a contradiction.

CASE 2.2) $\beta_{1}^{22}=0.$

Using the left-hand side of (\ref{form:bracketbeta}), we get $\beta_{2}^{22}=0.$

CASE 2.2.1) $\beta_{2}^{21}\neq0.$ In this case, by Proposition \ref{pro: normalization of beta},
we can assume $\beta_{2}^{21}=1.$

Using the left-hand side of (\ref{form:bracketbeta}), we get $c_{22}^{21}=0,c_{12}^{21}=0.$

Using the right-hand side of (\ref{form:bracketbeta}), we get a contradiction.

CASE 2.2.2) $\beta_{2}^{21}=0.$ Since $\dim\left(\mathrm{Im}\beta\right)=2$
we have $\beta_{2}^{12}\neq0.$ By Proposition \ref{pro: normalization of beta},
we can assume $\beta_{2}^{12}=1.$

Using (\ref{form:antisymmetrybeta}), we obtain $c_{12}^{21}=0,c_{12}^{12}=-1.$

Using the left-hand side of (\ref{form:bracketbeta}), we get $\beta_{1}^{21}=0$
which contradicts the condition $\dim\left(\mathrm{Im}\beta\right)=2.$ \end{proof}

\end{document}